\documentclass[10pt,twoside]{article}

\usepackage[T1]{fontenc}
\usepackage[english]{babel}

\usepackage{amsmath,amssymb,amsthm}
\usepackage{xcolor,textcomp,geometry,etoolbox}
\IfFileExists{newtxmath.sty}{\usepackage{newtxmath}}{%
  \usepackage{txfonts}

  \DeclareMathSizes{9}{9}{7.3}{5.5}
  \DeclareMathSizes{10}{10}{7.3}{5.5}
  \DeclareMathSizes{10.95}{10.95}{8}{6.1}
}
\usepackage[numbers]{natbib}

\makeatletter
\renewcommand\normalsize{%
  \@setfontsize\normalsize{10pt}{13pt}%
  \abovedisplayskip 6pt plus2pt minus0pt
  \abovedisplayshortskip 4pt plus3pt
  \belowdisplayshortskip 4pt plus3pt minus3pt
  \belowdisplayskip\abovedisplayskip
  \let\@listi\@listI}
\renewcommand\small{\@setfontsize\small{9pt}{11pt}}
\renewcommand\footnotesize{\@setfontsize\footnotesize{7pt}{8.5pt}}
\renewcommand\scriptsize{\@setfontsize\scriptsize{9pt}{9pt}}
\renewcommand\large{\@setfontsize\large{12pt}{14pt}}
\renewcommand\Large{\@setfontsize\Large{16pt}{18pt}}
\normalsize
\makeatother
\makeatletter
\def\paperheadfont{\fontsize{10bp}{13bp}\selectfont}
\def\ps@paper{%
  \let\@oddfoot\@empty\let\@evenfoot\@empty
  \def\@evenhead{\hbox to\textwidth{\rlap{\paperheadfont\thepage}\hfill
    {\paperheadfont\itshape J. Yu, S. Vaienti and W. Zhang}\hfill}}%
  \def\@oddhead{\hbox to\textwidth{\hfill
    {\paperheadfont\itshape Smooth rigidity for random expanding circle maps}\hfill
    \llap{\paperheadfont\thepage}}}}
\def\ps@plain{%
  \let\@oddfoot\@empty\let\@evenfoot\@empty
  \def\@oddhead{\hfill{\paperheadfont\thepage}}%
  \let\@evenhead\@oddhead}
\def\papersectionfont{\normalfont\fontsize{10bp}{13bp}\itshape\selectfont\raggedright}
\def\@seccntformat#1{\textup{\csname the#1\endcsname.}\hskip7.5pt}
\renewcommand\section{\@startsection{section}{1}{0pt}
  {-13bp plus-2pt minus2pt}{0.05bp}{\papersectionfont}}
\renewcommand\subsection{\@startsection{subsection}{2}{0pt}
  {-13pt plus-2pt minus2pt}{-1em}{\papersectionfont}}
\renewcommand\subsubsection{\@startsection{subsubsection}{3}{0pt}
  {-13pt plus-2pt minus2pt}{-1em}{\papersectionfont}}
\patchcmd{\@sect}{\@svsec #8}{\@svsec #8.}{}{
  \PackageError{manuscript}{Unable to format run-in headings}{}}
\renewcommand\@pnumwidth{1.55em}
\renewcommand\@tocrmarg{2.55em}
\renewcommand\@dotsep{4.5}
\renewcommand*\l@section[2]{%
  \ifnum\c@tocdepth>0
    \addpenalty\@secpenalty\addvspace{0pt}%
    \setlength\@tempdima{1.5em}%
    \begingroup
      \parindent0pt\rightskip\@pnumwidth\parfillskip-\@pnumwidth
      \leavevmode\advance\leftskip1.5em\hskip-\leftskip
      #1\nobreak\hfil\nobreak\hb@xt@\@pnumwidth{\hss #2}\par
    \endgroup
  \fi}

\renewenvironment{enumerate}{%
  \ifnum\@enumdepth>3\@toodeep\else
    \advance\@enumdepth\@ne
    \edef\@enumctr{enum\romannumeral\the\@enumdepth}%
    \list{\textup{(\arabic{\@enumctr})}}{%
      \usecounter{\@enumctr}\settowidth\labelwidth{0.}%
      \labelsep1em\leftmargin\labelwidth\advance\leftmargin\labelsep
      \topsep0pt\partopsep0pt\itemsep0pt\parsep0pt
      \listparindent\parindent\raggedright
      \def\makelabel##1{\hss\llap{##1}}}%
  \fi}{\endlist}

\newtheoremstyle{paperplain}{8.5pt}{8.5pt}{\normalfont\itshape}{0pt}
  {\scshape}{.}{.35em}{\thmname{#1}\thmnumber{ \textup{#2}}\thmnote{ {\normalfont #3}}}
\newtheoremstyle{paperremark}{8.5pt}{8.5pt}{\normalfont}{0pt}
  {\scshape}{.}{.35em}{\thmname{#1}\thmnumber{ \textup{#2}}\thmnote{ {\normalfont #3}}}
\newtheoremstyle{paperdefinition}{8.5pt}{8.5pt}{\normalfont}{0pt}
  {\normalfont}{.}{.35em}{}

\newif\ifpaperproof
\newtheoremstyle{paperproof}{8.5pt}{8.5pt}{\normalfont}{0pt}
  {\itshape}{.}{.45em}{\global\paperprooftrue
    \@ifempty{#3}{\thmname{#1}}{\thmname{#3}}}
\def\@endtheorem{\ifpaperproof\qed\global\paperprooffalse\fi
  \endtrivlist\@endpefalse}

\newcommand\paperabstracthead{\@startsection{section}{1}{0pt}
  {-17pt}{-5pt}{\normalfont\fontsize{10bp}{13bp}\itshape\selectfont}}
\renewenvironment{abstract}{\par\addvspace{0pt}
  \paperabstracthead*{Abstract.}\fontsize{10bp}{13bp}\selectfont}
  {\par\addvspace{13pt}}
\newcommand\keywords[2][Key words]{\par\noindent
  {\fontsize{10pt}{13pt}\selectfont #1: #2\par}}
\newenvironment{smallprint}{\par\fontsize{8pt}{10pt}\selectfont
  \setlength{\labelsep}{.5em}%
  \renewcommand\section{\@startsection{section}{1}{0pt}{-10.25pt}{2.5pt}
    {\normalfont\fontsize{10bp}{13bp}\scshape\selectfont\centering}}%
  }{\par}
\newenvironment{ack}{\subsection*{Acknowledgements.}\par}{\par}

\renewcommand\@biblabel[1]{[#1]\hfill}
\makeatother

\theoremstyle{paperplain}
\newtheorem{lm}{Lemma}
\newtheorem{thm}{Theorem}
\newtheorem{pro}{Proposition}
\newtheorem{cor}{Corollary}
\theoremstyle{paperdefinition}
\newtheorem{defi}{Definition}
\theoremstyle{paperremark}
\newtheorem{rmk}{Remark}
\theoremstyle{paperproof}
\newtheorem*{proof}{Proof}
\numberwithin{equation}{section}
\allowdisplaybreaks[4]

\usepackage{hyperref}
\hypersetup{colorlinks,breaklinks=true,plainpages=false,
  citecolor=black,linkcolor=black,urlcolor=black,
  bookmarksopen=true,bookmarksnumbered=false,bookmarksdepth=3,
  pdftitle={Smooth rigidity for randomly perturbed expanding maps on the circle},
  pdfauthor={Jiangqiong Yu, Sandro Vaienti, Wenmeng Zhang}}

\begin{document}
\raggedbottom
\thispagestyle{plain}

% Title and affiliations use ordinary LaTeX commands.
\noindent\vbox{\vskip51.5pt\hsize\textwidth\parindent0pt
  \hypersetup{urlcolor=blue}
  {\centering\fontsize{17.25bp}{22pt}\bfseries\selectfont\mathversion{bold}
    Rigidity of randomly perturbed expanding\\circle maps\par}
  \vskip15pt
  {\centering\fontsize{10bp}{13bp}\selectfont
    JIANGQIONG YU\textdagger, SANDRO VAIENTI\textdaggerdbl\ and WENMENG ZHANG\textdagger\par}
  \vskip3pt
  {\centering\fontsize{10bp}{13bp}\itshape\selectfont
    \textup{\textdagger}\ School of Mathematical Sciences, Chongqing Normal University,\\
    Chongqing 401331, P.~R.~China\\
    e-mail: \textcolor{blue}{\href{mailto:jiangqiongyu168@163.com}{jiangqiongyu168@163.com}} (Jiangqiong Yu)\\
    e-mail: \textcolor{blue}{\href{mailto:wmzhang@cqnu.edu.cn}{wmzhang@cqnu.edu.cn}} (Wenmeng Zhang)\par
    \addvspace{.5\baselineskip}
    \textup{\textdaggerdbl}\ Aix Marseille Universit\'e, Universit\'e de Toulon, CNRS, CPT,\\
    Marseille 13009, France\\
    and Great Bay University, Dongguan 523808, P.~R.~China\\
    e-mail: \textcolor{blue}{\href{mailto:sandro.vaienti@gmail.com}{sandro.vaienti@gmail.com}} (Sandro Vaienti)\par
    \addvspace{.5\baselineskip}
    {\normalfont\small Corresponding author: Wenmeng Zhang}\par
    \addvspace{.5\baselineskip}}
}\par\addvspace{43pt}

\begin{abstract}
We extend the classical smooth linearization theory of expanding circle maps to randomly perturbed ones without requiring small perturbations, establishing two equivalent rigidity criteria of absolute continuity of random conjugacy and Lyapunov exponents along random periodic orbits. In the random setting, an essential difficulty is the uniqueness of absolutely continuous random invariant measures.
For the first criterion, we overcome this difficulty by refining the memory-loss argument of Ott, Stenlund and Young, which gives the exponential convergence of probability densities, to establish the uniqueness of random invariant measures; for the second criterion, we bypass this difficulty by employing periodic-orbital averaging together with Jensen's equality condition to derive pointwise rigidity directly from the periodic data.
%without invoking uniqueness.

\end{abstract}



% Contents are printed in black.
\begingroup
\color{black}
\hypersetup{hidelinks}
\renewcommand{\color}[2][]{}
\renewcommand{\textcolor}[3][]{#3}
\let\bfseries\mdseries
\tableofcontents
\endgroup

\section{Introduction}
Smooth linearization is a long-standing classical problem in the study of dynamical systems, which is concerned with reducing a nonlinear system to a linear one via a smooth change of coordinates. As early as the 1950s-1960s, in the vicinity of a hyperbolic fixed point, Grobman and Hartman independently (\cite{Grobman,Hartman}) proved a local topological linearization theorem, while Sternberg (\cite{Sternberg}) established smooth linearization theorem under additional non-resonant conditions. For circle maps, the investigation can be traced back to the earlier 1930s, when Denjoy (\cite{Denjoy}) provided a topological conjugacy between an orientation-preserving circle diffeomorphism with irrational rotation number and a rigid rotation. The subsequent linearization results of Arnold (\cite{Arnold1961}), Herman (\cite{herman1979}) and Yoccoz (\cite{Yoccoz1984}) showed how arithmetic properties of the rotation number govern smooth conjugacy.
Those developments illustrate the essential distinctions between obtaining topological coordinates and establishing their higher regularity.

From the 1980s, the problem of establishing smoothness of the conjugacy from some weaker assumptions has become
a fundamental and widely studied problem, which is usually called the {\it rigidity} problem. One of the most basic framework for studying this problem is given by expanding maps on the circle $\mathbb{S}^1$. In 1985, based on the topological conjugacy result between two expanding circle maps (\cite{Shub1969}), Shub and Sullivan (\cite{Shub1985}) proved that if the two $C^r$ ($r\ge 2$) expanding circle maps admit an absolutely continuous conjugacy, then the prescribed topological conjugacy is $C^r$. Later in 1994,
Arteaga (\cite{Arteaga1994}) proved that if a \(C^{r,\gamma}\) ($0<\gamma<1$) expanding circle map satisfies the periodic-point derivative condition, then its topological conjugacy to the linear model
$E_d(x)=dx \pmod{1}$ is \(C^r\).
These results identify absolute continuity and periodic data as important criteria for rigidity.

On the $2$-dimensional torus $\mathbb T^2$, in the 1980s-1990s, de la Llave and Moriy{\'o}n (\cite{Llave1987, Llave-Mor1988, Llave}) published a series of works establishing rigidity result for smooth Anosov diffeomorphisms, where the smooth conjugacy was obtained from the existence of absolute continuous conjugacy or from the equality of Lyapunov exponents along periodic orbits. Higher-dimensional results include spectral rigidity for Anosov diffeomorphisms preserving conformal structures on their stable and unstable bundles (\cite{KaSa2009,Llave2002}), as well as for irreducible and non-algebraic systems (\cite{Gogolev2017,Gogolev,KaSaWa2023,SaYa2019}).

On the other hand, the study of expanding circle maps has been extended to the random setting, since this framework captures stochastic stability under small random perturbations and rich statistical behaviours such as mixing, limit theorems and equilibrium states (\cite{BaladiYoung1993, Kifer1986, Kifer1998, KhaninKifer1996}). In 2017, Nakano (\cite{Nakano2017}) proved a random version of Shub's topological conjugacy theorem for random expanding circle maps with small perturbations, which depends  continuously on the base point. Recently, Brown and Shi (\cite[Theorem~1.1]{BrownShi2026}) proved that the random walk generated by a finite family of $C^r$ expanding maps on $\mathbb{S}^1$ admits a single smooth diffeomorphism simultaneously conjugating every map to an affine map when the Lyapunov exponent of its unique stationary SRB-measure attains the average logarithm of the degrees.

In this paper, we investigate random circle maps (with sample point $\omega$) of a fixed degree
$|d|\ge 2$ on every $\omega$-fiber and their smooth conjugacy to the linear model
$E_{d}(x)$, where the fixed degree is a sufficient and necessary condition for the existence of topological random  conjugacy with the conjugacy equation holding on every fiber (see Proposition 1 below). These random circle maps can be regarded as perturbations of $E_{d}$ since they share the same degree ${d}$ as $E_{d}$, but no smallness of the perturbations is required. Our main results provide the following two
equivalent criteria characterizing when the prescribed topological random conjugacy becomes smooth.

\par\nobreak\smallskip\nobreak
(i) The existence of an
absolutely continuous random conjugacy almost
surely with respect to $\omega$, which will be proved in the following Theorem 1 (smooth rigidity);

\par\nobreak\smallskip\nobreak
(ii) the Lyapunov exponent at every random periodic orbit equals $\log d$ almost
surely with respect to $\omega$, which will be proved in the following Theorem 2 (spectral rigidity).
\par\nobreak\smallskip\nobreak

Note that the above smoothness concerns the state variable and that a random periodic orbit means a finite cyclic collection of graphs under the skew product sense (see Definition~\ref{def:random-periodic-graphs} below). Moreover,
the above results extend the classical rigidity theory of Shub-Sullivan (\cite{Shub1985}) and Arteaga (\cite{Arteaga1994}) to random circle dynamics, respectively. Yet, our proofs differ from their deterministic arguments as the uniqueness of the absolutely continuous invariant measure, a standard ingredient of the deterministic argument in \cite{Arteaga1994,Shub1985}, does not transfer automatically to the random setting and needs to be re-established.

To prove criterion (i), we refine the memory-loss argument of Ott, Stenlund and Young (\cite{OttStenlundYoung2009}), which gives the exponential convergence of probability densities, by combining it with introducing a \(C^1\) regularity invariant cone, and with using methods of smooth approximation of densities and smooth normalization of systems. This refining approach addresses the problems of uniqueness and smoothness for associated random invariant measures. Uniqueness then identifies the random invariant measure induced by the absolute continuous conjugacy with the one used for the smooth normalization, giving a Lebesgue measure preserving conjugacy between the normalized system and $E_d$. This forces the normalized system to be affine, which
together with the smooth normalizing transformation yields smoothness of the prescribed topological random conjugacy.

To prove criterion (ii), at every random periodic point, we average the logarithmic derivative over the base space, to reduce the random periodic Lyapunov condition to deterministic periodic data, and then use the Riemann sums over the dense set of periodic points to obtain an integral identity on $\mathbb{S}^1$.
Jensen's equality condition then turns this averaged constraint (i.e., the integral identity on $\mathbb{S}^1$) into pointwise rigidity, forcing the normalized system to have constant derivative, i.e., to be affine. It enables us to reduce the remainder of the proof to case (i).
This sufficiency argument bypasses the issue of uniqueness of invariant measures used in \cite{Arteaga1994}.

The paper is organized as follows. Section~\ref{sec:main} introduces the preliminaries and states the two main theorems and their corollary. Section~\ref{sec:rigidity-proof} states the density proposition and applies it to establish the two rigidity criteria. Section~\ref{sec:densities} gives a self-contained proof of the density proposition, covering its construction, uniqueness, and full regularity in the state variable. Section~\ref{sec:examples} presents examples illustrating the two criteria. Appendix~\ref{app:topological-linearization} contains the proof of the topological random linearization theorem, Appendix~\ref{app:periodic-graph-density} shows the density of random periodic points, and Appendix~\ref{app:higher-order-density} proves the higher-order derivative and H\"older estimates for the random density.

%%%%%%%%%%%%%%%%%%%%%%%%%%%%%%%%%%%%%%%%%%%%%%%%%%%%%%%%%%%%%%%%%%%%%%%%%%%%%%%
\section{Preliminaries and main results}\label{sec:main}

\subsection{Preliminaries}
Let $(\Omega,\mathcal F,\mathbb P)$ be a probability space, where $\mathcal F$ denotes the $\sigma$-algebra of measurable subsets of $\Omega$, and $\mathbb P$ is the probability measure defined on {$\mathcal F$}. Let $\theta:\Omega\to\Omega$ be a measure-preserving transformation, and
let $\mathbb S^1:=\mathbb R/\mathbb Z$, equipped with the usual distance
$d(x,y)$, the Borel \(\sigma\)-algebra $\mathcal B(\mathbb S^1)$ and normalized Lebesgue measure $m$.
%{For a metric space $X$, $\mathcal B(X)$ denotes its Borel $\sigma$-algebra.
%The abbreviations a.s. and a.e. mean almost surely and almost everywhere, respectively.}
We call $\psi_\omega:\mathbb S^1\to\mathbb S^1$ a
\emph{continuous random circle map} if
$\omega\mapsto\psi_\omega(x)$ is measurable for every
$x\in\mathbb S^1$ and $\psi_\omega$ is continuous for every
$\omega\in\Omega$.
These conditions imply joint measurability in $(\omega,x)$,
since $\mathbb S^1$ is a separable metric space.
Here, $\Omega$ is the sample space,
$\omega$ is the sample point, $x$ is the state variable and \((\Omega,\mathcal F,\mathbb P,\theta)\) is the base space.

For a continuous random circle map $\psi_\omega$, write
\begin{align}\label{cocycle-Psi}
\Psi(n,\omega)
:=\psi_{\theta^{n-1}\omega}\circ\cdots\circ\psi_\omega,
\quad \forall n\ge1,
\qquad \Psi(0,\omega)=\mathrm{Id},
\end{align}
where $\mathrm{Id}$ denotes the identity map.
We also write $\Psi(n,\omega,x):=\Psi(n,\omega)(x)$.
%when it exists, $D\Psi(n,\omega,x)$ denotes the derivative with respect to $x$.
These iterates satisfy the cocycle identity
\[
\Psi(n+k,\omega)
=\Psi(n,\theta^k\omega)\circ\Psi(k,\omega),
\qquad \forall n,k\ge0.
\]
For every fixed $n\ge0$, $\Psi(n,\omega)(x)$ is jointly
measurable in $(\omega,x)$ and continuous in $x$ for every $\omega$.
The copy of $\mathbb S^1$ over $\omega$ is a fiber, and
$\psi_\omega$ can be viewed as a fiber map.
For background on random dynamical systems, see \cite{Arnold1998}.

We call $h_\omega:\mathbb S^1\to\mathbb S^1$ a
\emph{random circle homeomorphism} if $h_\omega$ is a
homeomorphism for every $\omega\in\Omega$ and both
$\omega\mapsto h_\omega(x)$ and
$\omega\mapsto h_\omega^{-1}(x)$ are measurable for every
fixed $x\in\mathbb S^1$.
It is orientation preserving, respectively orientation reversing,
if every $h_\omega$ has that orientation.
Let $\psi_\omega,\tilde\psi_\omega:\mathbb S^1\to\mathbb S^1$
be continuous random circle maps over the same base space \((\Omega,\mathcal F,\mathbb P,\theta)\).
A random circle homeomorphism $h_\omega$ is called a
\emph{topological random conjugacy} from $\psi_\omega$ to
$\tilde\psi_\omega$ if
\[
h_{\theta\omega}\circ\psi_\omega
=\tilde\psi_\omega\circ h_\omega,
\quad \mathbb{P}\text{-a.s.}
\]
Here $\psi_\omega$ is the original random circle map and
$\tilde\psi_\omega$ is the target random circle map.

For an integer $K\ge1$, a continuous random circle map
$\psi_\omega:\mathbb S^1\to\mathbb S^1$ is called a
\emph{random $C^K$ circle map} if its derivatives $D^j\psi_\omega$,
$1\le j\le K$, exist and are continuous in the state variable for every
$\omega\in\Omega$.
For each $x$, the maps $\omega\mapsto D^j\psi_\omega(x)$ are
measurable, since derivatives are obtained successively as pointwise
limits of measurable difference quotients.
For $0<\alpha\le1$, a random $C^K$ circle map is called a
\emph{random $C^{K,\alpha}$ circle map} if $D^K\psi_\omega$ is
$\alpha$-H\"older continuous for every $\omega\in\Omega$.
It is called a \emph{random $C^\infty$ circle map} if it is a
random $C^K$ circle map for every integer $K\ge1$.

Fix an integer $d\ge2$ and define
\[
E_d(x)=dx\pmod1,
\qquad E_{-d}(x)=-dx\pmod1.
\]
We study the linearization problem in which the target map is $\tilde\psi_\omega=E_d$ or $\tilde\psi_\omega=E_{-d}$. In what follows,
the main statements and proofs are formulated for orientation-preserving
fiber maps with target map $E_d$, as the discussion for an orientation-reversing version is similar.

Assume that the random $C^1$ expanding circle map $\psi_\omega$
satisfies the following uniform expansion bound
\begin{equation}\label{eq:standing-expansion}
\inf_{\omega\in\Omega}\inf_{x\in\mathbb S^1}|D\psi_\omega(x)|\ge\lambda>1,
\end{equation}
where the derivative of $\psi_\omega$ is defined as the derivative of its lift $ \phi_\omega$.The lift $\phi_\omega:\mathbb R\to\mathbb R$ satisfies
$\phi_\omega(t)\bmod1=\psi_\omega(t\bmod1)$.
Orientation preservation then gives $D\psi_\omega(x)\ge\lambda$ for every $\omega\in\Omega$ and every $x\in\mathbb S^1$.
Let $\deg\psi_{\omega}$ be the degree of $\psi_\omega$ defined by its lift, i.e.,
$\deg\psi_{\omega}:= \phi_\omega(t+1)- \phi_\omega(t)$,
which is independent of $t\in\mathbb R$.
Besides, we write $D\Psi(n,\omega,x)$ denotes the derivative with respect to $x$.

Then, we have the topological random linearization result.

\begin{pro}[(Topological linearization)]\label{th:continuous-linearization}
Suppose that $\psi_\omega:\mathbb S^1\to\mathbb S^1$ is an
orientation-preserving random $C^1$ expanding circle map satisfying {\rm (\ref{eq:standing-expansion})}.
Then there exists an orientation-preserving random circle
homeomorphism $h^0_\omega:\mathbb S^1\to\mathbb S^1$ satisfying
\begin{equation}\label{eq:topological-conjugacy}
h^0_{\theta\omega}\circ\psi_\omega
=E_d\circ h^0_\omega,
\qquad \forall\omega\in\Omega,
\end{equation}
if and only if $\deg\psi_\omega=d$ for every $\omega\in\Omega$.
\end{pro}

The proof is a straightforward adaptation of the classical construction
using lifts to $\mathbb R$ and a contraction argument
(\cite[Section~I, Theorems~2--3]{Shub1969}) to the random setting.
For each fixed $x$, the successive approximations for the conjugacy and
its inverse are measurable in $\omega$, so their pointwise limits are
also measurable.
We include the proof in Appendix~\ref{app:topological-linearization}
for completeness. Moreover,
by Proposition~\ref{th:continuous-linearization}, it is natural to build our subsequent study of smooth conjugacy upon the framework of random systems with an identical degree $d$ on every $\omega$-fiber. Thus, the random circle maps under consideration can be regarded as perturbations of $E_{d}$, since they share the same degree ${d}$ as $E_{d}$, but no smallness of the perturbations is required.

\subsection{Smooth rigidity theorem}
Next,
we provide the first criterion concerning absolutely continuous random conjugacy to characterize when $h^0_\omega$ becomes smooth. For a real-valued function $v$ on $X=\mathbb S^1$ or $\mathbb R$ and
$0<\alpha\le1$, set
\[
\|v\|_{\infty}:=\sup_{x\in X}|v(x)|,
\qquad
[v]_{\alpha}:=
\sup_{\substack{x,y\in X\\x\ne y}}
\frac{|v(x)-v(y)|}{d_X(x,y)^\alpha},
\]
where $d_{\mathbb S^1}(x,y)=\min_{k\in\mathbb Z}|x-y-k|$ and
$d_{\mathbb R}(x,y)=|x-y|$. Assume that $\psi_\omega$ is a
random $C^{K,\alpha}$ circle map with $K\ge2$ and $0<\alpha\le1$,
and impose the uniform bounds
\begin{equation}\label{eq:regularity-hypotheses}
U:=\sup_{\omega\in\Omega}
\max_{2\le j\le K}\|D^j\psi_\omega\|_\infty<\infty,
\qquad
M:=\sup_{\omega\in\Omega}
[D^K\psi_\omega]_\alpha<\infty{.}
\end{equation}
The finite constants $U,M\ge0$ need not be small.

A finite Borel measure $\nu$ on $\mathbb S^1$ is
absolutely continuous with respect to $m$, written $\nu\ll m$,
if $m(A)=0$ implies $\nu(A)=0$ for every Borel set $A$.
We write $\nu\sim m$ if $\nu\ll m$ and $m\ll\nu$.

We call $h_\omega$ an \emph{absolutely continuous random conjugacy}
from $\psi_\omega$ to $\tilde\psi_\omega$ if the following
conditions hold:

\par\nobreak\smallskip\nobreak
(A1) For $\mathbb{P}$-almost every $\omega$, $h_\omega$ is a
Borel isomorphism between two Borel subsets of $\mathbb S^1$
of full $m$-measure, and both $h_\omega$ and $h_\omega^{-1}$
send $m$-null sets to $m$-null ones such that
\begin{equation}\label{eq:ac-conjugacy}
h_{\theta\omega}(\psi_\omega(x))
=\tilde\psi_\omega(h_\omega(x))
\quad\text{for }m\text{-almost every }x\in\mathbb S^1.
\end{equation}

(A2) The pushforward measure
\begin{equation}\label{eq:transported-family}
\nu_\omega=(h_\omega^{-1})_*m,\quad \mathbb P\text{-a.s.},
\end{equation}
where $((h_\omega^{-1})_*m)(A):=m(h_\omega(A))$ for $A\in\mathcal B(\mathbb S^1)$, satisfies that $\omega\mapsto\nu_\omega(A)$ is measurable.
\par\nobreak\smallskip\nobreak

Condition (A1) extends the measure-theoretic notion of
absolutely continuous conjugacy used by Shub and Sullivan
\cite{Shub1985} to the random setting.
No continuity or orientation condition is imposed on $h_\omega$.
For the Borel isomorphisms in (A1), preservation of $m$-null
sets in both directions is equivalent to
$\nu_\omega=(h_\omega^{-1})_*m\sim m$
for $\mathbb P$-almost every $\omega$.
Since joint measurability of the map
$(\omega,x)\mapsto h_\omega^{-1}(x)$ is not assumed,
(A2) explicitly requires that
$\omega\mapsto\nu_\omega(A)$ be measurable for every
$A\in\mathcal B(\mathbb S^1)$.
This ensures that $\nu_\omega$ is a random probability measure and can be integrated over the base.

\iffalse
The measurability requirement in~\eqref{eq:transported-family} specifies the required dependence of the transported probabilities on $\omega$.
It is automatic if the inverse family $(h_\omega^{-1})$ has a jointly measurable representative $G:\Omega\times {\mathbb{S}^1}\to {\mathbb{S}^1}$: for almost every $\omega$,
$G(\omega,\cdot)=h_\omega^{-1}$ almost everywhere, and
\[
\nu_\omega(B)=\int_{{\mathbb{S}^1}}\mathbf1_B(G(\omega,y))\,dm(y),
\qquad B\in\mathcal B({\mathbb{S}^1}),
\]
defines a measurable probability family on all of $\Omega$.
In particular, for a family of circle homeomorphisms whose inverses are measurable in the base point $\omega$ for each fixed $x\in {\mathbb{S}^1}$, continuity in the state variable implies joint measurability and hence this requirement.
The initial maps $h_\omega$ above, however, are only Borel isomorphisms modulo null sets; their continuity in the state variable and joint measurability are not assumed. We therefore retain the measurable-probability condition rather than impose joint measurability on representatives of these maps.
\fi

In what follows, we always let \((\Omega,\mathcal F,\mathbb P,\theta)\) be the base space given at the beginning of this section with invertible $\theta$ such that $\theta$ and $\theta^{-1}$ are measurable. Then,
our first main theorem concerning the absolute-continuity criterion is the following.

\begin{thm}[(Smooth rigidity)]\label{th2.2}
Let $d\ge2$, $K\ge2$ be integers and $0<\alpha\le1$.
Suppose that $\psi_\omega:\mathbb{S}^1\to \mathbb{S}^1$ is an orientation-preserving random $C^{K,\alpha}$ expanding circle map, with $\deg\psi_\omega=d$ for every $\omega\in\Omega$.
Assume that
\[
\begin{gathered}
\inf_{\omega\in\Omega}\inf_{x\in\mathbb S^1}
|D\psi_\omega(x)|\ge\lambda>1,~~
\sup_{\omega\in\Omega}
\left(\max_{2\le j\le K}\|D^j\psi_\omega\|_\infty
      +[D^K\psi_\omega]_\alpha\right)<\infty.
\end{gathered}
\]
Let $h^0_\omega:\mathbb{S}^1\to \mathbb{S}^1$ be the topological random conjugacy obtained by Proposition~{\rm \ref{th:continuous-linearization}}.
Then, $h^0_\omega$ is an orientation-preserving random $C^{K,\alpha}$ diffeomorphism for $\mathbb P$-almost every $\omega$ if and only if $\psi_\omega $ admits an absolutely continuous random conjugacy to $E_d$. Moreover,
if each $\psi_\omega$ is $C^\infty$ and
\begin{equation}\label{eq:cinfty-uniform-bounds}
\sup_{\omega\in\Omega}\|D^j\psi_\omega\|_\infty<\infty
\qquad\text{for every integer }j\ge2,
\end{equation}
the same equivalence holds with $C^\infty$ in place of $C^{K,\alpha}$.
\end{thm}

In general, the absolutely continuous conjugacy need not coincide with $h^0_\omega$ (see e.g. \cite{Shub1985}). The proof of Theorem~\ref{th2.2} will be given in Section~\ref{sec:rigidity-proof}.

\subsection{Spectral rigidity theorem}
Then, we provide the second criterion concerning Lyapunov exponents along random periodic orbits to characterize when $h^0_\omega$ becomes smooth.
Consider the skew product
\[
T:\Omega\times\mathbb{S}^1\longrightarrow\Omega\times\mathbb{S}^1,
\qquad
T(\omega,x)=(\theta\omega,\psi_\omega(x)).
\]
For a measurable section $\xi:\Omega\to\mathbb{S}^1$, we write
$
\Gamma_\xi:=\{(\omega,\xi(\omega)):\omega\in\Omega\}.
$

Following~\cite[Definition~2.2]{HLL}, we describe random
periodic orbits in terms of measurable graphs that are
cyclically permuted by the skew product $T$.

\begin{defi}\label{def:random-periodic-graphs}
Let $p\ge1$ and let
$Y_0,\ldots,Y_{p-1}:\Omega\to\mathbb S^1$ be measurable maps.
The collection of graphs $\{\Gamma_{Y_i}:0\le i<p\}$
is called a \emph{random periodic orbit} of $T$ with
period $p$ if
\begin{equation}\label{eq:random-periodic-graphs}
T\Gamma_{Y_i}=\Gamma_{Y_{(i+1)\bmod p}},
\quad \forall 0\le i<p,
\end{equation}
or, equivalently,
$$
\widetilde T(Y_i)=Y_{(i+1)\bmod p},\quad \forall 0\le i<p,
$$
where $\widetilde T$ is the induced map on measurable sections given by
$
\widetilde T(g)(\omega)
:=\psi_{\theta^{-1}\omega}(g(\theta^{-1}\omega)).
$
Each graph $\Gamma_{Y_i}$ is called a
\emph{random periodic point} of $T$; if no confusion arises, we also call the map $Y_i$ a \emph{random periodic point}.
\end{defi}

\begin{rmk}\label{rmk:periodic-graph-density}
In \cite[Theorem~2.1(i)]{HLL}, for a quasi-periodically
forced Anosov diffeomorphism, it is proved that its random
periodic points are dense in the space of measurable sections
with an appropriate distance.
The same conclusion also holds
for an orientation-preserving random $C^1$ circle map
satisfying~\eqref{eq:standing-expansion}
in the space $\mathcal M(\Omega,\mathbb S^1)$ of measurable
sections with the distance
\begin{align}\label{dis-section}
d_\infty(g,\xi):=
\sup_{\omega\in\Omega}d(g(\omega),\xi(\omega)),
\end{align}
which will be proved in Appendix~\ref{app:periodic-graph-density}.
In particular, the values $Y_i(\omega)$ of random periodic points
are dense in the $\omega$-fiber.
\end{rmk}

For a random periodic orbit $Y$ of period $p$, define its
fiber Lyapunov exponents by
\[
\lambda_{Y_i}(\omega):=
\lim_{n\to\infty}\frac{1}{n}
\log D\Psi(n,\omega,Y_i(\omega)),
\qquad 0\le i<p.
\]
Under the hypotheses below, these limits exist
$\mathbb{P}$-almost surely for each fixed $Y$ and all
its phases, by Birkhoff's ergodic theorem applied
to $\theta^p$ and the uniform derivative bounds.
Let $\psi_\omega:\mathbb S^1\to\mathbb S^1$
satisfy the same hypotheses as the part before ``Assume that" in Theorem \ref{th2.2}. Then,
our second main theorem concerning the Lyapunov exponent criterion is the following.

\iffalse
For a random periodic orbit $Y$ of period $p$, its fiber
Lyapunov exponents are defined, whenever the limits exist, by
\[
\lambda_{Y_i}(\omega):=
\lim_{n\to\infty}\frac1n
\log D\Psi(n,\omega,Y_i(\omega)),
\qquad 0\le i<p,
\]
where $\Psi$ is defined in (\ref{cocycle-Psi}). Under the hypotheses below, the derivative cocycle over the
measure-preserving phase extension
$S_Y(\omega,i)=(\theta\omega,i+1\bmod p)$,
with uniform probability on the phases, satisfies the
integrability conditions of the multiplicative ergodic
theorem (\cite{Arnold1998}).
Thus, for each fixed $Y$, these limits exist
$\mathbb P$-almost surely, simultaneously for all phases.
\fi

\begin{thm}[(Spectral rigidity)]
\label{th:lyapunov-bootstrap}
Assume that
\[
\begin{gathered}
\inf_{\omega\in\Omega}
\inf_{x\in\mathbb S^1}|D\psi_\omega(x)|
\ge\lambda>1,~~
\sup_{\omega\in\Omega}
\left(
\max_{2\le j\le K}\|D^j\psi_\omega\|_\infty
+[D^K\psi_\omega]_\alpha
\right)<\infty.
\end{gathered}
\]
Let $h^0_\omega$ be the topological random conjugacy
obtained by Proposition~\ref{th:continuous-linearization}.
Then $h^0_\omega$ is an orientation-preserving random
$C^{K,\alpha}$ diffeomorphism for $\mathbb P$-almost every
$\omega$ if and only if
\begin{equation}\label{eq:lyapunov-criterion}
\lambda_{Y_i}(\omega):=
\lim_{n\to\infty}\frac1n
\log D\Psi(n,\omega,Y_i(\omega))=\log d
\end{equation}
for every random periodic orbit $Y$ of period $p$
in Definition~{\rm \ref{def:random-periodic-graphs}},
every $0\le i<p$, and $\mathbb P$-almost every $\omega$.
Moreover, if $\psi_\omega$ is a random $C^\infty$ circle map and
$\sup_{\omega\in\Omega}\|D^j\psi_\omega\|_\infty<\infty$, for every integer $j\ge2$, the same equivalence holds with $C^\infty$
in place of $C^{K,\alpha}$.
\end{thm}

The proof of Theorem~\ref{th:lyapunov-bootstrap} will be given in Section~\ref{sec:rigidity-proof}.
The equivalence of the above two criteria is given in the following corollary.

\iffalse
\begin{rmk}[Uniform bounds in the random setting]\label{rmk:uniform-bounds}
For a fixed $C^{K,\alpha}$ circle map, the relevant derivative norms and H\"older seminorm in the state variable are finite on the compact circle.
For a random family, finiteness for each fixed base point does not imply a bound uniform over all $\omega\in\Omega$.
The uniform expansion bound controls inverse branches, while the bounds on higher derivatives with respect to the state variable and on H\"older seminorms provide uniform distortion and regularity estimates in the state variable for the densities used in the smoothness criteria.
They are sufficient hypotheses for the present proofs; no smallness of $U$ or $M$ is assumed, and no higher derivative bounds are needed for Proposition~\ref{th:continuous-linearization}.
In the notation~\eqref{eq:lift-notation}, $D^jF_\omega(t)=D^j\psi_\omega(t\bmod1)$ for $2\le j\le K$.
A bound $\sup_\omega\|DF_\omega\|_\infty<d-1$ would imply uniform expansion, but is stronger than the direct assumption~\eqref{eq:standing-expansion} and is not imposed here.
No bound on the original constants $F_\omega(0)$ is needed: integer normalization of a lift supplies a uniformly bounded displacement, as shown in the proof of Proposition~\ref{th:continuous-linearization}.
\end{rmk}
\fi

\begin{cor}\label{cor:ac-lyapunov}
Suppose that $\psi_\omega:\mathbb{S}^1\to \mathbb{S}^1$ is an orientation-preserving random $C^{K,\alpha}$ expanding circle map, with $\deg\psi_\omega=d$ for every $\omega\in\Omega$.
Then, $\psi_\omega$ admits an absolutely continuous random
conjugacy to $E_d$ if and only if the periodic Lyapunov
condition~\eqref{eq:lyapunov-criterion} holds for every random
periodic orbit for $\mathbb P$-almost every $\omega$.
\end{cor}

\begin{proof}
By Theorems~\ref{th2.2} and~\ref{th:lyapunov-bootstrap},
each condition is equivalent to $C^{K,\alpha}$ regularity
of the same prescribed topological conjugacy $h^0_\omega$
for $\mathbb P$-almost every $\omega$.
\end{proof}

\section{Proofs of the two main theorems}\label{sec:rigidity-proof}

The sufficiency proofs of both theorems follow a common strategy.
We first use a random invariant density to construct a smooth change of coordinates
that sends the random invariant measure to Lebesgue measure.
In these coordinates, the hypothesis of each theorem will imply that the
fiber maps have constant derivative $d$.
An affine rigidity argument then shows that the prescribed topological conjugacy
differs from the density coordinates only by rotations, and hence has the required regularity.

To carry out this strategy, we need the invariant-density
proposition and the affine rigidity lemma stated below.
Before presenting these results, we introduce the necessary
measure-theoretic notation
(see~\cite{Arnold1998,Crauel2002}).

We call $\mu_\omega$ a \emph{random invariant measure} (or an \emph{equivariant random probability measure}) on $\mathbb S^1$ if
\begin{enumerate}
\item[(1)] for every $A\in\mathcal B(\mathbb S^1)$, the map $\omega\mapsto\mu_\omega(A)\in[0,1]$ is $\mathcal F$-measurable;
\item[(2)] for $\mathbb P$-almost every $\omega$, the map $A\mapsto\mu_\omega(A)$ is a Borel probability measure on $\mathbb S^1$;
\item[(3)] $(\psi_\omega)_*\mu_\omega=\mu_{\theta\omega}$ for $\mathbb P$-almost every $\omega$.
\end{enumerate}
Here $(\psi_\omega)_*\mu_\omega(A):=\mu_\omega(\psi_\omega^{-1}(A))$ for $A\in\mathcal B(\mathbb S^1)$. Conditions (1)--(2) define a \emph{random probability measure}. See~\cite[Theorem~1.4.5(ii)]{Arnold1998}.

A nonnegative jointly measurable function
$(\omega,x)\mapsto\rho_\omega(x)$ with
$\int_{\mathbb S^1}\rho_\omega\,dm=1$ for $\mathbb P$-almost every $\omega$
defines a random probability measure $\mu_\omega:=\rho_\omega m$;
we call $\rho_\omega$ a \emph{random probability density}.
%We also write $\rho(\omega,x):=\rho_\omega(x)$ for this jointly measurable function.
Here $\mu_\omega=\rho_\omega m$ means that
\[
\mu_\omega(A)=\int_A\rho_\omega\,dm,
\qquad A\in\mathcal B(\mathbb S^1).
\]
For $\mathbb P$-almost every $\omega$, $\rho_\omega$ is the
Radon--Nikodym derivative $d\mu_\omega/dm$,
uniquely determined up to $m$-null sets.
$L^1(\mathbb S^1,m)$ denotes the space of $m$-integrable functions modulo equality $m$-almost everywhere.
For an orientation-preserving random $C^1$ expanding circle map
$\psi_\omega$, its transfer operator {(see~\cite{StenlundSulku2014})} acts on $v\in L^1(\mathbb S^1,m)$ by
\begin{equation}\label{eq:transfer-operator}
(\mathcal L_\omega v)(x)
:=\sum_{\psi_\omega(z)=x}
\frac{v(z)}{D\psi_\omega(z)},
\qquad m\text{-a.e. }x\in\mathbb S^1.
\end{equation}
We define
\[
\mathcal L_\omega^n
:=\mathcal L_{\theta^{n-1}\omega}\circ\cdots\circ\mathcal L_\omega,
\quad \forall n\ge1,
\qquad \mathcal L_\omega^0=\mathrm{Id}.
\]
The identity
$
(\psi_\omega)_*(\rho_\omega m)
=(\mathcal L_\omega\rho_\omega)m
$
shows that $\mu_\omega=\rho_\omega m$ is invariant if and only if
$\mathcal L_\omega\rho_\omega=\rho_{\theta\omega}$ in
$L^1(\mathbb S^1,m)$ for $\mathbb P$-almost every $\omega$.

%For measures $\mu$ and $\nu$, the notation $\mu\ll\nu$ means that $\mu$ is absolutely continuous with respect to $\nu$, and $\mu\sim\nu$ means that each is absolutely continuous with respect to the other.

Next, for an integer $r\ge0$ and $0<\alpha\le1$, we use the norm
\[
\|v\|_{C^{r,\alpha}}
:=\sum_{j=0}^{r}\|D^jv\|_\infty+[D^rv]_\alpha,
\qquad D^0v:=v.
\]
The following proposition supplies the positive random invariant density and the
uniform regularity bounds needed for this change of coordinates.

\begin{pro}\label{pro1}
Suppose that $\psi_\omega:\mathbb S^1\to\mathbb S^1$
is an orientation-preserving random $C^{K,\alpha}$
{\rm (}$K\ge2$, $0<\alpha\le1${\rm )} expanding map with $\deg\psi_\omega=d\ge2$ for every $\omega\in\Omega$.
Assume that
\[
\begin{gathered}
\inf_{\omega\in\Omega}
\inf_{x\in\mathbb S^1}|D\psi_\omega(x)|
\ge\lambda>1,~~
\sup_{\omega\in\Omega}
\left(
\max_{2\le j\le K}\|D^j\psi_\omega\|_\infty
+[D^K\psi_\omega]_\alpha
\right)<\infty.
\end{gathered}
\]
Then there exist constants $c_0,C_0>0$ and a random invariant probability density
$\rho_\omega:\mathbb S^1\to(0,\infty)$ such that
$\rho_\omega\in C^{K-1,\alpha}(\mathbb S^1)$,
$\int_{\mathbb S^1}\rho_\omega\,dm=1$,
and $c_0\le\rho_\omega\le C_0$ for every $\omega\in\Omega$. Moreover,
\begin{equation}\label{eq:uniform-density-regularity}
\sup_{\omega\in\Omega}
\|\rho_\omega\|_{C^{K-1,\alpha}}<\infty.
\end{equation}
The density $\rho_\omega$ defines a random invariant
probability measure $\mu_\omega:=\rho_\omega m$.
If $\bar\mu_\omega$ is another random invariant probability
measure with $\bar\mu_\omega\ll m$ for
$\mathbb P$-almost every $\omega$, then
$\bar\mu_\omega=\mu_\omega$ for
$\mathbb P$-almost every $\omega$.
\end{pro}

The proof of Proposition~\ref{pro1} will be given in
Section~\ref{sec:densities}.
In the sufficiency proofs of both theorems, we use the following lemma
to show that an orientation-preserving topological random conjugacy
from an affine system to $E_d$ is a rotation on almost every fiber.
%Once the normalized maps are affine, the following lemma identifies their orientation-preserving conjugacies to $E_d$ as rotations.

\begin{lm}\label{lm:affine-rigidity}
Let $T_\omega(x)=dx+a_\omega\pmod1$, where $d\ge2$
and $a_\omega\in\mathbb R$.
For each $\omega$, let $q_\omega$ be an orientation-preserving
circle homeomorphism, and suppose that
\[
q_{\theta\omega}\circ T_\omega=E_d\circ q_\omega
\]
for $\mathbb P$-almost every $\omega$.
Then $q_\omega$ is a rotation for
$\mathbb P$-almost every $\omega$.
\end{lm}

\begin{proof}
Since $\theta$ preserves $\mathbb P$, for almost every $\omega$
the conjugacy relation holds at $\theta^j\omega$ for every $j\ge0$.
Fix such an $\omega$ and choose an increasing degree-one lift
$\widehat q_{\theta^j\omega}$ of $q_{\theta^j\omega}$
for each $j\ge0$.
For $n\ge1$, composing the real affine lifts $t\mapsto dt+a_{\theta^j\omega}$, $j=0,\ldots,n-1$, gives
\[
t\longmapsto d^nt+a_{\omega,n},
\qquad
a_{\omega,n}
:=\sum_{j=0}^{n-1}d^{\,n-1-j}a_{\theta^j\omega}.
\]
Here $a_{\omega,n}$ is the accumulated translation over $n$ steps, from the fiber over $\omega$ to the fiber over $\theta^n\omega$.
Iterating the conjugacy relation shows that
\[
\widehat q_{\theta^n\omega}(d^nt+a_{\omega,n})
-d^n\widehat q_\omega(t)
\]
is an integer independent of $t$.
Subtracting its values at $t$ and $s$ gives
\begin{align*}
&d^n\bigl[(\widehat q_\omega(t)-t)
          -(\widehat q_\omega(s)-s)\bigr]
\\
&=\bigl[\widehat q_{\theta^n\omega}(d^nt+a_{\omega,n})
                  -(d^nt+a_{\omega,n})\bigr]
-\bigl[\widehat q_{\theta^n\omega}(d^ns+a_{\omega,n})
                  -(d^ns+a_{\omega,n})\bigr].
\end{align*}
For any increasing degree-one lift $\widehat q$, the periodic function $\widehat q-\mathrm{Id}$ has oscillation at most one.
Indeed, after translating one argument by an integer, we may assume $s\le t\le s+1$, so both $\widehat q(t)-\widehat q(s)$ and $t-s$ belong to $[0,1]$.
Consequently,
\[
\left|(\widehat q_\omega(t)-t)-(\widehat q_\omega(s)-s)\right|
\le d^{-n}.
\]
Letting $n\to\infty$ shows that $\widehat q_\omega-\mathrm{Id}$ is constant.
Thus $q_\omega$ is a rotation.
This completes the proof of Lemma 1.
\end{proof}

\subsection{Proof of the smooth rigidity theorem}
\begin{proof}[Proof of Theorem~\ref{th2.2}]
For necessity, suppose that $h^0_\omega$ is an orientation-preserving $C^{K,\alpha}$ diffeomorphism
for almost every $\omega$. Then $h^0_\omega$ and its inverse preserve $m$-null sets.
By Proposition~\ref{th:continuous-linearization}, the inverse map $(h^0_\omega)^{-1}$ is
measurable in $\omega$ at each fixed point and continuous in the state variable,
hence jointly measurable. Thus, for every $A\in\mathcal B(\mathbb S^1)$,
\[
\nu^0_\omega(A)
:=\bigl(((h^0_\omega)^{-1})_*m\bigr)(A)
=\int_{{\mathbb{S}^1}}\mathbf1_A\bigl((h^0_\omega)^{-1}(y)\bigr)\,dm(y)
\]
is measurable in $\omega$. Here $\mathbf1_A$ denotes the indicator function of $A$.
Together with~\eqref{eq:topological-conjugacy}, these facts verify (A1) and (A2).

For sufficiency, let $h_\omega$ be an absolutely continuous random conjugacy to $E_d$,
and let $\mu_\omega=\rho_\omega m$ be supplied by Proposition~\ref{pro1}.

\par\medskip
{\itshape Step 1. Identify the random invariant measure obtained by $h_\omega$ with
the one obtained in Proposition \ref{pro1}, which is used to construct smooth normalizing coordinates.}
\par\nobreak\smallskip\nobreak

By (A1) and (A2), $\nu_\omega=(h_\omega^{-1})_*m$
is a random probability measure that is absolutely continuous
with respect to $m$ for $\mathbb P$-almost every $\omega$.
The conjugacy relation and $(E_d)_*m=m$ imply
\[
(\psi_\omega)_*\nu_\omega
=(h_{\theta\omega}^{-1}\circ E_d)_*m
=(h_{\theta\omega}^{-1})_*m
=\nu_{\theta\omega}
\]
for $\mathbb P$-almost every $\omega$, with compositions
understood modulo $m$-null sets.
Uniqueness in Proposition~\ref{pro1} identifies the transported measure with $\mu_\omega$:
\begin{equation}\label{eq:identified-measure}
(h_\omega^{-1})_*m=\mu_\omega,
\qquad
(h_\omega)_*\mu_\omega=m
\end{equation}
for $\mathbb P$-almost every $\omega$.

To construct a smooth change of coordinates sending $\mu_\omega$ to $m$,
extend $\rho_\omega$ periodically to $\mathbb R$ and set
\begin{equation}\label{eq:density-coordinate}
\widehat h^1_\omega(t)=\int_0^t \rho_\omega(s)\,ds,
\qquad
h^1_\omega(t\bmod1)
=\widehat h^1_\omega(t)\bmod1.
\end{equation}
Since $\int_{\mathbb S^1}\rho_\omega\,dm=1$ and $\rho_\omega>0$, we have
\[
\widehat h^1_\omega(t+1)
=\widehat h^1_\omega(t)+1,
\qquad
D\widehat h^1_\omega(t)=\rho_\omega(t)>0.
\]
Thus $h^1_\omega$ is an orientation-preserving
$C^{K,\alpha}$ circle diffeomorphism with
$h^1_\omega(0)=0$.
A change of variables gives
\begin{equation}\label{eq:normalization}
(h^1_\omega)_*\mu_\omega=m,
\qquad
\bigl((h^1_\omega)^{-1}\bigr)_*m=\mu_\omega.
\end{equation}

Joint measurability of $\rho$ makes
$\omega\mapsto\widehat h^1_\omega(t)$ measurable for each fixed $t$.
The inverse lifts are measurable at each fixed $y$, since, for every $r\in\mathbb Q$,
\[
(\widehat h^1_\omega)^{-1}(y)<r
\quad\Longleftrightarrow\quad
y<\widehat h^1_\omega(r).
\]
Continuity in the state variable then gives joint measurability
of $(\omega,x)\mapsto h^1_\omega(x)$ and $(\omega,y)\mapsto(h^1_\omega)^{-1}(y)$.
Define
\begin{equation}\label{eq:normalized-maps}
f_\omega
=h^1_{\theta\omega}\circ\psi_\omega
 \circ(h^1_\omega)^{-1}.
\end{equation}
Thus $f_\omega$ is a random $C^{K,\alpha}$ circle map with positive
derivative and degree $d$ for every $\omega\in\Omega$.

The invariance of the random measure $\mu_\omega$
and~\eqref{eq:normalization} give, for
$\mathbb P$-almost every $\omega$,
\[
(f_\omega)_*m
=(h^1_{\theta\omega})_*(\psi_\omega)_*\mu_\omega
=(h^1_{\theta\omega})_*\mu_{\theta\omega}
=m.
\]
The inverse-branch formula and continuity therefore yield
\[
\sum_{f_\omega(z)=y}\frac{1}{Df_\omega(z)}=1,
\qquad y\in\mathbb S^1.
\]
Since $d\ge2$ and all summands are positive, each is strictly
less than one. Hence $Df_\omega(x)>1$ for every $x$, and
compactness gives $\min_{x\in\mathbb S^1}Df_\omega(x)>1$.
Thus $f_\omega$ is an orientation-preserving random
$C^{K,\alpha}$ expanding circle map preserving Lebesgue
measure, for $\mathbb P$-almost every $\omega$.

\par\medskip
\noindent{\itshape Step 2. Prove that the normalized maps have constant derivative.}
\par\nobreak\smallskip\nobreak

For almost every $\omega$, define
\[
h^2_\omega=h_\omega\circ(h^1_\omega)^{-1}.
\]
By (A1) and the smoothness of $h^1_\omega$, this is a Borel
isomorphism between conull Borel sets.
The conjugacy relations give
\begin{equation}\label{eq:normalized-conjugacy}
h^2_{\theta\omega}\circ f_\omega
=E_d\circ h^2_\omega
\quad m\text{-almost everywhere}.
\end{equation}
Equations~\eqref{eq:identified-measure}
and~\eqref{eq:normalization} give
\[
\begin{aligned}
(h^2_\omega)_*m
=(h_\omega)_*
  \bigl(\bigl((h^1_\omega)^{-1}\bigr)_*m\bigr)
  =(h_\omega)_*\mu_\omega
=m.
\end{aligned}
\]
Since $\theta$ preserves $\mathbb P$, for almost every $\omega$
the conjugacy relation~\eqref{eq:normalized-conjugacy} holds and both $h^2_\omega$ and $h^2_{\theta\omega}$ are
measure-preserving Borel isomorphisms. Fix such an $\omega$.
There is a conull Borel set $E\subset {\mathbb{S}^1}$
on which $h^2_\omega$ is defined,
$f_\omega(x)$ belongs to the domain of
$h^2_{\theta\omega}$, and
\eqref{eq:normalized-conjugacy} holds.
Indeed, the inverse image under $f_\omega$
of the exceptional domain set has measure zero,
since $(f_\omega)_*m=m$.

To determine the derivative of $f_\omega$, we compare the measures of sets on which $f_\omega$ is injective.
Let $I$ be an arc on which $f_\omega$ is injective,
and let $A\subset I\cap E$ be Borel.
If $x_1,x_2\in A$ satisfy
\[
E_d(h^2_\omega(x_1))=E_d(h^2_\omega(x_2)),
\]
then \eqref{eq:normalized-conjugacy} gives
\[
h^2_{\theta\omega}(f_\omega(x_1))
=h^2_{\theta\omega}(f_\omega(x_2)).
\]
Injectivity of $h^2_{\theta\omega}$ and then of
$f_\omega$ on $I$ imply $x_1=x_2$.
Hence $E_d$ is injective on $h^2_\omega(A)$.

For any Borel set $B$ on which $E_d$ is injective,
\[
m(E_d(B))=d\,m(B).
\]
Indeed, partition $B$ by the $d$ half-open intervals
$[j/d,(j+1)/d)$, $0\le j<d$.
On each interval the map multiplies length by $d$,
and the images of these pieces are disjoint
by injectivity.

Applying this fact to $B=h^2_\omega(A)$ and using
the measure-preserving properties of $h^2_\omega$
and $h^2_{\theta\omega}$, we obtain
\begin{equation}\label{eq:jacobian-sets}
\begin{aligned}
m(f_\omega(A))
=m\bigl(h^2_{\theta\omega}(f_\omega(A))\bigr)
=m\bigl(E_d(h^2_\omega(A))\bigr)
=d\,m(h^2_\omega(A))
=d\,m(A).
\end{aligned}
\end{equation}
The maps $h^2_\omega$ and $h^2_{\theta\omega}$ are Borel isomorphisms, while $f_\omega|_I$ and $E_d|_{h^2_\omega(A)}$ are continuous and injective. Hence the images in~\eqref{eq:jacobian-sets} are Borel.

For an arbitrary Borel set $A\subset I$, replace
$A$ by $A\cap E$.
The removed set has measure zero, and its image
under the $C^1$ map $f_\omega$ also has measure zero.
Hence \eqref{eq:jacobian-sets} holds for every
Borel $A\subset I$.

Since $f_\omega$ is injective on $I$ and has
positive derivative, change of variables yields
\[
\int_A Df_\omega\,dm
=m(f_\omega(A))
=d\,m(A),
\qquad A\in\mathcal B(I).
\]
It follows that $Df_\omega=d$ almost everywhere on $I$ and,
by continuity, everywhere. Covering the circle by such arcs gives
\begin{equation}\label{eq:constant-jacobian}
Df_\omega(x)=d,
\qquad x\in {\mathbb{S}^1}.
\end{equation}
This is the constant-Jacobian argument of \cite{Shub1985}, expressed through the measures of injective sets.
No derivative of $h_\omega$ or $h^2_\omega$ with respect to the state variable is used.

\par\medskip
\noindent{\itshape Step 3. Prove the smoothness of the prescribed topological random conjugacy.}
\par\nobreak\smallskip\nobreak

We now compare the prescribed topological conjugacy $h^0_\omega$ with the smooth density coordinate $h^1_\omega$.
Equation~\eqref{eq:constant-jacobian} shows that
$f_\omega$ is affine for almost every $\omega$. The maps
$h^0_\omega\circ(h^1_\omega)^{-1}$
are orientation-preserving circle homeomorphisms.
Their conjugacy relations imply
\[
\begin{aligned}
\bigl(h^0_{\theta\omega}
      \circ(h^1_{\theta\omega})^{-1}\bigr)
      \circ f_\omega
&=
E_d\circ
\bigl(h^0_\omega\circ(h^1_\omega)^{-1}\bigr).
\end{aligned}
\]

Apply Lemma~\ref{lm:affine-rigidity} to the affine maps
$x\mapsto dx+f_\omega(0)\pmod1$, which agree with $f_\omega$ almost surely.
It follows that $h^0_\omega\circ(h^1_\omega)^{-1}$ is a rotation
for $\mathbb P$-almost every $\omega$.
Since $h^1_\omega(0)=0$, the rotation angle is $h^0_\omega(0)$ as a circle element.
Consequently,
\begin{equation}\label{eq:bootstrap-identification}
h^0_\omega(x)
=h^1_\omega(x)+h^0_\omega(0)\pmod1,
\qquad
Dh^0_\omega(x)=\rho_\omega(x).
\end{equation}
Thus the previously constructed $h^0_\omega$ is an orientation-preserving $C^{K,\alpha}$ diffeomorphism for $\mathbb P$-almost every $\omega$.
For the $C^\infty$ assertion, the uniform bound on $D^{K+1}\psi_\omega$ controls the Lipschitz seminorm of $D^K\psi_\omega$ at each finite $K\ge2$.
Apply the finite-order result with $\alpha=1$ for every $K$ and intersect the resulting countably many sets of full probability for the same conjugacy $h^0_\omega$.
The converse follows because a smooth circle diffeomorphism and its inverse preserve null sets.
This completes the proof of Theorem 1.
\end{proof}

\begin{rmk}
Absolute continuity ensures that $\nu_\omega=(h_\omega^{-1})_*m\ll m$
almost surely, hence the existence of $d\nu_\omega/dm$;
continuity of $h_\omega$ alone does not guarantee $\nu_\omega\ll m$.
Without this absolute continuity, the uniqueness assertion in Proposition~\ref{pro1} cannot be
applied to obtain $\nu_\omega=\rho_\omega m$. This identity is needed to make
$h_\omega\circ(h^1_\omega)^{-1}$ measure preserving.
\end{rmk}

\subsection{Proof of the spectral rigidity theorem}
\begin{proof}[Proof of Theorem~\ref{th:lyapunov-bootstrap}]

Let $\rho_\omega$ be the density in Proposition~\ref{pro1},
with $c_0\le\rho_\omega\le C_0$, and set $\nu^0_\omega=((h^0_\omega)^{-1})_*m$.
Joint measurability of $(\omega,x)\mapsto(h^0_\omega)^{-1}(x)$
makes $\nu^0_\omega$ a random probability measure.
The conjugacy equation and $(E_d)_*m=m$ give
{
\begin{equation}\label{eq:topological-measure-invariance}
(\psi_\omega)_*\nu^0_\omega=\nu^0_{\theta\omega}.
\end{equation}}
Moreover, $\nu^0_\omega$ has full support, since
$\nu^0_\omega(V)=m(h^0_\omega(V))>0$
for every nonempty open set $V\subset\mathbb S^1$.

For necessity, suppose that $h^0_\omega$ is an orientation-preserving $C^{K,\alpha}$ diffeomorphism for almost every $\omega$.
On these fibers, change of variables gives $\nu^0_\omega=(Dh^0_\omega)m$.
Uniqueness in Proposition~\ref{pro1} and continuity give
$Dh^0_\omega(x)=\rho_\omega(x)$ for every $x\in\mathbb S^1$,
for $\mathbb P$-almost every $\omega$.
Since $\theta$ preserves $\mathbb P$, for almost every $\omega$
this equality and the stated smoothness hold at $\theta^n\omega$
for every $n\ge0$. Fix such an $\omega$.
Differentiating the conjugacy equation gives
\begin{equation}\label{eq:density-derivative-criterion}
\rho_{\theta\omega}(\psi_\omega(x))D\psi_\omega(x)=d\rho_\omega(x),
\qquad x\in {\mathbb{S}^1}.
\end{equation}
Iteration yields
\begin{equation}\label{eq:iterated-density-criterion}
D\Psi(n,\omega,x)
=d^n\frac{\rho_\omega(x)}
 {\rho_{\theta^n\omega}(\Psi(n,\omega,x))}.
\end{equation}
Therefore, for $\mathbb P$-almost every $\omega$, the following estimates hold for all $x\in {\mathbb{S}^1}$ and $n\ge1$:
\begin{equation}\label{eq:uniform-orbit-growth}
\begin{gathered}
\frac{c_0}{C_0}d^n\le D\Psi(n,\omega,x)
                  \le\frac{C_0}{c_0}d^n,
                  \quad
\left|\frac1n\log D\Psi(n,\omega,x)-\log d\right|
\le\frac{\log(C_0/c_0)}n.
\end{gathered}
\end{equation}
In particular, substituting $x=Y_i(\omega)$ gives \eqref{eq:lyapunov-criterion} for every random periodic orbit and each of its phases, proving necessity.

For sufficiency, suppose that~\eqref{eq:lyapunov-criterion} holds.

\par\medskip
\noindent{\itshape Step 1. Derive an averaged logarithmic derivative identity from periodic data.}
\par\nobreak\smallskip\nobreak

For $p\ge1$, write $\operatorname{Fix}(E_d^p)=\{a\in\mathbb S^1:E_d^p(a)=a\}$ for the fixed-point set, and set
\[
\mathcal P_d=\bigcup_{p\ge1}\operatorname{Fix}(E_d^p),
\qquad
\xi_a(\omega)=(h^0_\omega)^{-1}(a),
\quad a\in\mathcal P_d.
\]
Each $\xi_a$ is measurable.
The conjugacy equation gives
\begin{equation}\label{eq:periodic-section-covariance}
\Psi(j,\omega,\xi_a(\omega))
=\xi_{E_d^j a}(\theta^j\omega),
\qquad j\ge0,\quad\omega\in\Omega.
\end{equation}
For $E_d^p a=a$, the sections
\begin{equation}\label{eq:induced-random-periodic-section-family}
Y_i^a(\omega)=\xi_{E_d^i a}(\omega),
\qquad 0\le i<p,
\end{equation}
form a random periodic orbit of period $p$.
The assumed condition thus applies to $\xi_a$ for every $a\in\mathcal P_d$.
%Averaging over the base converts the random periodic condition into deterministic periodic data, from which uniform Riemann sums yield an integral identity (cf.~\cite[Section~2, Remark~3]{OHare2025}).

For $y\in {\mathbb{S}^1}$, put
\[
g_\omega(y)=\log D\psi_\omega((h^0_\omega)^{-1}(y)),
\qquad
\overline g(y)=\int_\Omega g_\omega(y)\,d\mathbb P(\omega).
\]
The function $g_\omega$ is continuous for every $\omega$,
and $(\omega,y)\mapsto g_\omega(y)$ is jointly measurable.
Since $\int_{\mathbb S^1}D\psi_\omega\,dm=d$ and
$\|D^2\psi_\omega\|_\infty\le U$, the mean value theorem
and~\eqref{eq:standing-expansion} give
\begin{equation}\label{eq:uniform-first-derivative-bounds}
\log\lambda\le g_\omega(y)\le\log(d+U).
\end{equation}
Dominated convergence therefore shows that $\overline g$ is continuous.
Fix $p\ge1$ and $a\in\operatorname{Fix}(E_d^p)$.
The covariance relation~\eqref{eq:periodic-section-covariance} and the chain rule give
\[
\log D\Psi(kp,\omega,\xi_a(\omega))
=\sum_{j=0}^{kp-1}g_{\theta^j\omega}(E_d^j a).
\]
Since $\theta$ preserves $\mathbb P$ and $E_d^p a=a$, integration gives
\[
\frac1{kp}\int_\Omega
 \log D\Psi(kp,\omega,\xi_a(\omega))\,d\mathbb P(\omega)
=\frac1p\sum_{j=0}^{p-1}\overline g(E_d^j a).
\]
By~\eqref{eq:uniform-first-derivative-bounds}, the normalized logarithm
is uniformly bounded, while the assumed criterion gives its almost-sure
limit $\log d$. Dominated convergence as $k\to\infty$ yields
\begin{equation}\label{eq:periodic-average-data}
\frac1p\sum_{j=0}^{p-1}\overline g(E_d^j a)=\log d,
\qquad a\in\operatorname{Fix}(E_d^p).
\end{equation}

For each $p\ge1$,
\[
\operatorname{Fix}(E_d^p)
=\left\{\frac{j}{d^p-1}:0\le j<d^p-1\right\}.
\]
The map $E_d$ permutes this finite set. Averaging \eqref{eq:periodic-average-data} over it therefore gives
\[
\frac1{d^p-1}
\sum_{a\in\operatorname{Fix}(E_d^p)}\overline g(a)=\log d.
\]
These are uniform Riemann sums for the continuous function $\overline g$. Letting $p\to\infty$ and using Fubini's theorem
and the definition of $\nu^0_\omega$, we obtain
\begin{equation}\label{eq:integrated-growth}
\begin{aligned}
\log d
=\int_{{\mathbb{S}^1}}\overline g(y)\,dm(y)
=\int_\Omega\int_{{\mathbb{S}^1}}
 \log D\psi_\omega(x)\,d\nu^0_\omega(x)\,d\mathbb P(\omega).
\end{aligned}
\end{equation}

\par\medskip
\noindent{\itshape Step 2. Prove that the normalized maps have derivative $d$ using Jensen's equality condition.}
\par\nobreak\smallskip\nobreak

%We use Jensen's equality condition to deduce that the normalized maps have derivative $d$ from the averaged identity~\eqref{eq:integrated-growth} (cf.~\cite[Proposition~28]{ContrerasLopesThieullen2001}).
Let $h^1_\omega$ be the density coordinate in~\eqref{eq:density-coordinate},
whose measurability and $C^{K,\alpha}$ regularity were established
in the proof of Theorem~\ref{th2.2}, and set
\[
f_\omega=h^1_{\theta\omega}\circ\psi_\omega\circ(h^1_\omega)^{-1}
\]
as in~\eqref{eq:normalized-maps}. The chain rule gives
\[
Df_\omega(h^1_\omega(x))
=\frac{\rho_{\theta\omega}(\psi_\omega(x))D\psi_\omega(x)}
 {\rho_\omega(x)}.
\]

For every continuous $v:\mathbb S^1\to\mathbb R$,
the uniform inverse-branch formula for $E_d$ and the conjugacy give
\begin{align}\label{eq:balanced-topological-measure}
\int_{{\mathbb{S}^1}}v(z)\,d\nu^0_\omega(z)
&=\int_{{\mathbb{S}^1}}v((h^0_\omega)^{-1}(\bar{z}))\,dm(\bar{z})
=\int_{{\mathbb{S}^1}}\frac1d
  \sum_{E_d\bar{z}=\bar{y}}v((h^0_\omega)^{-1}(\bar{z}))\,dm(\bar{y})
  \nonumber\\
&=\int_{{\mathbb{S}^1}}\frac1d
  \sum_{\psi_\omega(z)=y}v(z)\,d\nu^0_{\theta\omega}(y).
\end{align}
For the last equality, put $y=(h^0_{\theta\omega})^{-1}(\bar{y})$;
the preimages correspond under $z=(h^0_\omega)^{-1}(\bar{z})$.

For $\psi_\omega(z)=y$, define the positive weights
\[
w_\omega(z;y)
=\frac{\rho_\omega(z)}{\rho_{\theta\omega}(y)D\psi_\omega(z)}
=\frac{1}{Df_\omega(h^1_\omega(z))}.
\]
For almost every $\omega$, the density equivariance $\mathcal L_\omega\rho_\omega=\rho_{\theta\omega}$ and continuity give,
for every $y\in\mathbb S^1$,
\[
\sum_{\psi_\omega(z)=y}w_\omega(z;y)
=\frac{\mathcal L_\omega \rho_\omega(y)}{\rho_{\theta\omega}(y)}=1.
\]
On these fibers, Jensen's inequality for the strictly concave logarithm gives
\begin{equation}\label{eq:jensen-gap}
J_\omega(y):=-\frac1d\sum_{\psi_\omega(z)=y}
                  \log\bigl(d\,w_\omega(z;y)\bigr)\ge0.
\end{equation}
Equality holds exactly when all these weights equal $1/d$.
Local inverse branches show that $J_\omega$ is continuous;
the inverse-branch measurability argument in the proof of Lemma~\ref{lm4.1} then gives
joint measurability of $J(\omega,y):=J_\omega(y)$. The uniform bounds on $\rho_\omega$
and $D\psi_\omega$ make $J$ bounded.

To show that $J_\omega\equiv0$ for $\mathbb P$-almost every $\omega$,
we compute its average over the fibers and the base.
Apply~\eqref{eq:balanced-topological-measure} to the continuous function
\[
z\longmapsto
\log\frac{\rho_{\theta\omega}(\psi_\omega(z))D\psi_\omega(z)}
 {d\rho_\omega(z)}.
\]
Using the equivariance~\eqref{eq:topological-measure-invariance} in the logarithmic density terms gives
\[
\begin{aligned}
\int_{{\mathbb{S}^1}}J_\omega\,d\nu^0_{\theta\omega}
&=\int_{{\mathbb{S}^1}}\log D\psi_\omega\,d\nu^0_\omega-\log d+\int_{{\mathbb{S}^1}}\log \rho_{\theta\omega}\,d\nu^0_{\theta\omega}
       -\int_{{\mathbb{S}^1}}\log \rho_\omega\,d\nu^0_\omega.
\end{aligned}
\]
All four integrals are bounded and measurable in $\omega$.
Integrating over $\Omega$, the last two terms cancel since
$\theta_*\mathbb P=\mathbb P$.
Equation~\eqref{eq:integrated-growth} then gives
\[
\int_\Omega\int_{{\mathbb{S}^1}}J_\omega(y)\,
 d\nu^0_{\theta\omega}(y)\,d\mathbb P(\omega)=0.
\]
For almost every $\omega$, $J_\omega\ge0$ and
$\int J_\omega\,d\nu^0_{\theta\omega}=0$.
Continuity of $J_\omega$ and full support of $\nu^0_{\theta\omega}$
therefore imply $J_\omega(y)=0$ for every $y$ on these fibers.
The equality case in~\eqref{eq:jensen-gap} now gives $w_\omega(z;y)=1/d$ for every preimage $z$ of every $y$.
Since $w_\omega(z;y)=1/Df_\omega(h^1_\omega(z))$ and $h^1_\omega$ is onto, we obtain
\begin{equation}\label{eq:spectral-constant-derivative}
Df_\omega(x)=d,\qquad x\in\mathbb S^1,
\quad\text{for }\mathbb P\text{-almost every }\omega.
\end{equation}
Thus $f_\omega$ is affine on these fibers.

\par\medskip
\noindent{\itshape Step 3. Prove the smoothness of the prescribed topological random conjugacy.}
\par\nobreak\smallskip\nobreak

By~\eqref{eq:spectral-constant-derivative}, the maps $f_\omega$ are affine for almost every $\omega$.
The orientation-preserving homeomorphisms
$h^0_\omega\circ(h^1_\omega)^{-1}$ conjugate $f_\omega$ to $E_d$.
Apply Lemma~\ref{lm:affine-rigidity} to the affine maps
$x\mapsto dx+f_\omega(0)\pmod1$, which agree with $f_\omega$
almost surely. These homeomorphisms are therefore rotations
for $\mathbb P$-almost every $\omega$. Since $h^1_\omega(0)=0$,
\[
h^0_\omega(x)=h^1_\omega(x)+h^0_\omega(0)\pmod1,
\qquad Dh^0_\omega(x)=\rho_\omega(x).
\]
Consequently the same $h^0_\omega$ is an orientation-preserving $C^{K,\alpha}$ diffeomorphism for almost every $\omega$.

For the $C^\infty$ assertion, apply the finite-order equivalence with $\alpha=1$ at every $K\ge2$, using the bound on $D^{K+1}\psi_\omega$ for the required Lipschitz seminorm.
Intersecting the countably many full-probability sets makes this same $h^0_\omega$ smooth to every order.
Necessity follows already from the finite-order necessity argument.
This completes the proof of Theorem 2.
\end{proof}

\section{{Proof of the proposition for random density}}\label{sec:densities}

This section proves Proposition~\ref{pro1} under its stated
hypotheses on the base transformation and the fiber maps,
independently of the two bootstrap criteria.
The proof begins with a uniform memory-loss estimate in
the $C^1$ norm on a cone of positive probability densities.
This estimate yields an invariant density as a $C^1$ limit
of densities obtained by evolving the constant density
from the past to the current fiber.
Smooth approximation extends the memory-loss property
to integrable functions of zero mean.
Together with preservation of the base measure, this
yields uniqueness among absolutely continuous random
invariant probability measures. Finally, uniform derivative and H\"older estimates are
established for the approximating densities.
These estimates, together with compactness and the
$C^1$ convergence already obtained, yield
$C^{K-1,\alpha}$ regularity of the limiting invariant
density, with bounds uniform in the base point.

To prepare the estimates used in this argument, recall that
$\deg\psi_\omega=d$ for every $\omega\in\Omega$,
as assumed in Proposition~\ref{pro1}.
Consequently, $\deg\Psi(n,\theta^{-n}\omega)=d^n$,
so each $n$-step map has exactly $d^n$ local inverse branches
over any sufficiently short arc.
Set $\varepsilon:=\lambda^{-1}\in(0,1)$.
The constants in the uniform density bounds may be chosen as
\begin{equation}\label{eq:cone-constants}
A=1+\frac{2U\varepsilon^2}{1-\varepsilon},
\qquad
c_0=e^{-A/2},
\qquad
C_0=e^{A/2}.
\end{equation}

\subsection{{The memory-loss argument}}
Let $\mathbf1_\omega$ denote the constant density one on $\mathbb{S}^1$, attached to the base point $\omega$.
We use the $C^1$ norm in the state variable
$\|v\|_{C^1}:=\|v\|_\infty+\|Dv\|_\infty$.

Every fiber is represented by the same circle $\mathbb{S}^1$, equipped with normalized Lebesgue measure $m$.
With the constant $A$ fixed by~\eqref{eq:cone-constants}, set
\[
C_A=
\left\{\rho\in C^1({\mathbb{S}^1}):
\rho>0,\quad \int_{{\mathbb{S}^1}}\rho\,dm=1,\quad
|D\rho|\le A\rho\right\}.
\]
This is the normalized slice of the positive cone defined by $|D\rho|\le A\rho$.
Under the above identification, $C_A$ is the same set on every fiber: neither its defining conditions nor $A$ depends on $\omega$.
We therefore use the single notation $C_A$, without an argument for the base point.
We write $\rho_\omega$ for a function on $\mathbb{S}^1$ attached to the base point $\omega$; the subscript records only the fiber.
Membership $\rho_\omega\in C_A$ does not by itself impose random invariance or measurability in the base point $\omega$.
The constant density $\mathbf1_\omega$ belongs to $C_A$ on every fiber.

We now establish a uniform decay estimate for differences of densities by adapting the memory-loss argument of Ott, Stenlund and Young (\cite[Section~2.2]{OttStenlundYoung2009}).
We use the normalized $C^1$ cone $C_A$ to obtain decay in the $C^1$ norm.
The resulting estimate gives $C^1$ convergence of the pullback densities and provides the starting point for the higher-order regularity argument.

\begin{lm}\label{lm:memory-loss}
Under the standing hypotheses of this section, the transfer operators
preserve the normalized cone:
\begin{equation}\label{eq:cone-invariance}
\mathcal L_\omega C_A\subset C_A
\qquad(\omega\in\Omega).
\end{equation}
Set $\tau=c_0/2$, fix $N\ge1$ such that $\varepsilon^N\le1/3$,
and put $C_*=2(1+A)C_0$. For every $\omega\in{\Omega}$,
$n\ge0$ and $\rho_\omega,\widetilde\rho_\omega\in C_A$,
\begin{equation}\label{eq:cone-memory-loss}
\|\mathcal L_\omega^n(\rho_\omega-\widetilde\rho_\omega)\|_{C^1}
\le C_*(1-\tau)^{\lfloor n/N\rfloor}.
\end{equation}
Here $\lfloor t\rfloor$ denotes the greatest integer not exceeding $t$.
The constants are independent of the starting fiber and the two
initial densities.
\end{lm}

Inequality (\ref{eq:cone-memory-loss}) tells that the \(C^1\) difference between any two starting densities shrinks uniformly and exponentially as we apply the transfer operator repeatedly. The system gradually ``forgets'' its initial state, hence one uses the name ``memory-loss''.

\begin{proof}
For a real-valued integrable function $\rho_\omega$ on $\mathbb{S}^1$, the operator $\mathcal L_\omega$ produces a function on the same circle attached to the base point $\theta\omega$.
Changing variables on each of the ${d}$ inverse branches gives
\[
\int_{{\mathbb{S}^1}}\mathcal L_\omega \rho_\omega\,dm
=\int_{{\mathbb{S}^1}}\rho_\omega\,dm.
\]
For $v\in L^1(\mathbb S^1,m)$, we write $\|v\|_1:=\int_{\mathbb S^1}|v(x)|\,dm(x)$.
This norm is taken on a fixed fiber; there is no integration over the base point $\omega$.
Consequently,
\begin{equation}\label{eq:direct-L1-norm}
\|\mathcal L_\omega \rho_\omega\|_1\le\|\rho_\omega\|_1,
\qquad \rho_\omega\in L^1({\mathbb{S}^1},m).
\end{equation}

For $\rho_\omega\in C_A$, $|D\log\rho_\omega|\le A$. Integrating along a shortest arc yields
\[
|\log\rho_\omega(x)-\log\rho_\omega(y)|
\le A\,d(x,y)\le A/2.
\]
Thus $e^{-A/2}\rho_\omega(y)\le\rho_\omega(x) \le e^{A/2}\rho_\omega(y)$. Integrating in $y$ gives
\begin{equation}\label{eq:cone-bounds}
c_0\le\rho_\omega\le C_0,\qquad
\|D\rho_\omega\|_\infty\le AC_0,\qquad
\|\rho_\omega\|_{C^1}\le(1+A)C_0.
\end{equation}
In particular, $C_A$ is closed in $C^1(\mathbb S^1)$:
the common lower bound $c_0>0$, normalization and the derivative
inequality all pass to a $C^1$ limit.

We next check that $\mathcal L_\omega$ preserves $C_A$ while taking a density at $\omega$ to one at $\theta\omega$, both defined on $\mathbb{S}^1$.
For each local inverse branch $z=z(x)$, we have $Dz=1/D\psi_\omega(z)$. Differentiating the transfer-operator formula~\eqref{eq:transfer-operator} therefore gives
\begin{equation}\label{eq:transfer-derivative}
D(\mathcal L_\omega\rho_\omega)(x)
=\sum_{\psi_\omega(z)=x}
\left(
\frac{D\rho_\omega(z)}{D\psi_\omega(z)^2}
-\frac{\rho_\omega(z)D^2\psi_\omega(z)}{D\psi_\omega(z)^3}
\right).
\end{equation}
Since $D\psi_\omega\ge\varepsilon^{-1}$ and
$|D^2\psi_\omega|\le U$, it follows that
\[
\begin{aligned}
|D(\mathcal L_\omega\rho_\omega)(x)|
&\le\sum_{\psi_\omega(z)=x}
\left(\frac{|D\rho_\omega(z)|}{D\psi_\omega(z)^2}
+\frac{\rho_\omega(z)|D^2\psi_\omega(z)|}
 {D\psi_\omega(z)^3}\right)\\
&\le(\varepsilon A+U\varepsilon^2)
 (\mathcal L_\omega\rho_\omega)(x).
\end{aligned}
\]
The choice of $A$ gives
\[
\varepsilon A+U\varepsilon^2
=\varepsilon A+\frac{(1-\varepsilon)(A-1)}2
=A-\frac{(1-\varepsilon)(A+1)}2<A.
\]
Positivity and preservation of the integral therefore imply
$
\mathcal L_\omega C_A
\subset C_A.
$

The uniform positive lower bound permits removal of the same
mass from any two densities. Set $\tau=c_0/2=\tfrac12e^{-A/2}$ and fix $N\ge1$ such that $\varepsilon^N\le1/3$.
For $\rho_\omega\in C_A$, put
\begin{equation}\label{eq:remaining-density}
\widehat\rho_\omega
=\frac{\rho_\omega-\tau\mathbf1_\omega}{1-\tau}.
\end{equation}
This is a strictly positive probability density on the same fiber.
Since $\rho_\omega\ge2\tau$,
\[
\frac{|D\widehat\rho_\omega|}{\widehat\rho_\omega}
=\frac{|D\rho_\omega|}{\rho_\omega-\tau}
\le\frac{A\rho_\omega}{\rho_\omega-\tau}\le2A.
\]
We next show that $N$ iterates map the density $\widehat\rho_\omega$ defined in~\eqref{eq:remaining-density} into $C_A$.
For $N$ iterates, branchwise differentiation and the chain rule give
\[
\begin{aligned}
&D(\mathcal L_\omega^N\widehat\rho_\omega)(x)=\sum_{\Psi(N,\omega,z)=x}
\frac{\widehat\rho_\omega(z)}{D\Psi(N,\omega,z)}
\Biggl[
\frac{D\widehat\rho_\omega(z)}
 {\widehat\rho_\omega(z)D\Psi(N,\omega,z)}\\
&\hspace{25mm}
-\sum_{j=0}^{N-1}
\frac{D^2\psi_{\theta^j\omega}(\Psi(j,\omega,z))}
 {D\psi_{\theta^j\omega}(\Psi(j,\omega,z))}
\frac{D\Psi(j,\omega,z)}{D\Psi(N,\omega,z)}
\Biggr].
\end{aligned}
\]
Using $|D\widehat\rho_\omega|/\widehat\rho_\omega\le2A$,
$|D^2\psi_{\theta^j\omega}|/D\psi_{\theta^j\omega}\le U\varepsilon$,
and the expansion bounds
\[
\frac{1}{D\Psi(N,\omega,z)}\le\varepsilon^N,
\qquad
\frac{D\Psi(j,\omega,z)}{D\Psi(N,\omega,z)}
\le\varepsilon^{N-j}\quad(0\le j<N),
\]
we obtain, after reindexing the geometric sum and using
\eqref{eq:cone-constants} and $\varepsilon^N\le1/3$,
\[
\begin{aligned}
|D(\mathcal L_\omega^N\widehat\rho_\omega)|
&\le\left(2A\varepsilon^N
+U\varepsilon^2\sum_{j=0}^{N-1}\varepsilon^j\right)
 \mathcal L_\omega^N\widehat\rho_\omega\\
&=\left(\frac{A-1}{2}
+\frac{3A+1}{2}\varepsilon^N\right)
 \mathcal L_\omega^N\widehat\rho_\omega\le\left(A-\frac13\right)
 \mathcal L_\omega^N\widehat\rho_\omega.
\end{aligned}
\]
Thus, for every $\rho_\omega\in C_A$, the density defined in~\eqref{eq:remaining-density} satisfies
\begin{equation}\label{eq:cone-propagation}
\mathcal L_\omega^N
\left(\frac{\rho_\omega-\tau\mathbf1_\omega}{1-\tau}\right)
\in C_A,
\qquad \rho_\omega\in C_A.
\end{equation}
The constants are independent of the starting fiber, so this
statement applies to every block of $N$ consecutive iterates.

To compare $\rho_\omega,\widetilde\rho_\omega
\in C_A$, set
$\rho^{(0)}_\omega=\rho_\omega$ and
$\widetilde\rho^{(0)}_\omega=\widetilde\rho_\omega$.
Define the successive remaining densities by
\[
\rho^{(j+1)}_{\theta^{(j+1)N}\omega}
=\mathcal L_{\theta^{jN}\omega}^{N}
 \left(\frac{\rho^{(j)}_{\theta^{jN}\omega}
 -\tau\mathbf1_{\theta^{jN}\omega}}{1-\tau}\right),
\qquad j\ge0,
\]
and use the same definition for
$\widetilde\rho^{(j+1)}_{\theta^{(j+1)N}\omega}$.
By~\eqref{eq:cone-propagation}, both densities at stage $j$
are functions on $\mathbb{S}^1$ attached to $\theta^{jN}\omega$ and belong to the same set $C_A$.
At the first stage,
\[
\begin{aligned}
\mathcal L_\omega^N\rho_\omega
 &=\tau\mathcal L_\omega^N\mathbf1_\omega
 +(1-\tau)\rho^{(1)}_{\theta^N\omega},
 \quad
\mathcal L_\omega^N\widetilde\rho_\omega
 &=\tau\mathcal L_\omega^N\mathbf1_\omega
 +(1-\tau)\widetilde\rho^{(1)}_{\theta^N\omega}.
\end{aligned}
\]
The common terms cancel on subtraction. This uses linearity
and does not require $\mathcal L_\omega\mathbf1_\omega
=\mathbf1_{\theta\omega}$.
Repeating the subtraction proves, for every $q\ge0$,
\[
\mathcal L_\omega^{qN}(\rho_\omega-\widetilde\rho_\omega)
=(1-\tau)^q
 \left(\rho^{(q)}_{\theta^{qN}\omega}
 -\widetilde\rho^{(q)}_{\theta^{qN}\omega}\right).
\]
For $n=qN+\ell$, where $0\le\ell<N$, it follows that
\[
\begin{aligned}
\mathcal L_\omega^n(\rho_\omega-\widetilde\rho_\omega)
&=(1-\tau)^q\left(
 \mathcal L_{\theta^{qN}\omega}^{\ell}
       \rho^{(q)}_{\theta^{qN}\omega}
-\mathcal L_{\theta^{qN}\omega}^{\ell}
       \widetilde\rho^{(q)}_{\theta^{qN}\omega}\right).
\end{aligned}
\]
Each of the two densities in parentheses belongs to
$C_A$ by \eqref{eq:cone-invariance}, so
\eqref{eq:cone-bounds} and the triangle inequality give~\eqref{eq:cone-memory-loss}. This estimate is uniform in $\omega$,
$n\ge0$ and the two initial densities in $C_A$.
This completes the proof of Lemma 2.
\end{proof}

\subsection{{The uniqueness of the random density}}
The memory-loss estimate gives convergence of the pullback densities.
Smooth approximation and preservation of $\mathbb P$ then yield uniqueness
among absolutely continuous random invariant probability measures.

\begin{lm}\label{lm4.1}
The pullback densities
\[
\rho_{\omega,0}=\mathbf1_\omega,
\qquad
\rho_{\omega,n}
=\mathcal L_{\theta^{-n}\omega}^{\,n}\mathbf1_{\theta^{-n}\omega}
\quad(n\ge1)
\]
converge in the $C^1$ norm in the state variable, uniformly over the base points $\omega\in \Omega$, to a random probability density $\rho_\omega$ satisfying
\[
c_0\le \rho_{\omega,n},\rho_\omega\le C_0,
\qquad
\|D\rho_{\omega,n}\|_\infty,
\|D\rho_\omega\|_\infty\le AC_0.
\]
Writing $\mu_\omega:=\rho_\omega m$, we have $\mathcal L_\omega \rho_\omega=\rho_{\theta\omega}$
for every $\omega\in\Omega$.
If $\bar\mu_\omega$ is any random probability measure such that
\[
\bar\mu_\omega\ll m,
\qquad
(\psi_\omega)_*\bar\mu_\omega=\bar\mu_{\theta\omega}
\quad\text{for }\mathbb P\text{-almost every }\omega,
\]
then $\bar\mu_\omega=\mu_\omega$ for $\mathbb P$-almost every $\omega$.
\end{lm}

\begin{proof}
Since $\mathbf1_\omega\in C_A$, the cone invariance~\eqref{eq:cone-invariance} implies
$\rho_{\omega,n}\in C_A$ for every $n\ge0$.
The bounds in~\eqref{eq:cone-bounds} therefore apply to all these
pullback densities.

For the pullback sequence and $p\ge1$, the cocycle identity gives
\[
\rho_{\omega,n+p}-\rho_{\omega,n}
=\mathcal L_{\theta^{-n}\omega}^{\,n}
 \bigl(\rho_{\theta^{-n}\omega,p}-\mathbf1_{\theta^{-n}\omega}\bigr).
\]
Both densities in parentheses lie in $C_A$ on the fiber over
$\theta^{-n}\omega$. The memory-loss estimate~\eqref{eq:cone-memory-loss}
therefore gives a uniform Cauchy bound for the pullback sequence:
\begin{equation}\label{eq:direct-pullback-Cauchy}
\|\rho_{\omega,n+p}-\rho_{\omega,n}\|_{C^1}
\le C_*(1-\tau)^{\lfloor n/N\rfloor}.
\end{equation}
The right-hand side tends to zero independently of $p$ and $\omega$. Completeness of $C^1({\mathbb{S}^1})$ therefore gives a limit $\rho_\omega$, with convergence in the $C^1$ norm in the state variable uniform over the base points $\omega\in \Omega$.
Passing to the limit preserves the integral, the bounds
$c_0\le \rho_\omega\le C_0$, and
$|D\rho_\omega|\le A\rho_\omega$.
Thus $\rho_\omega\in C_A$.
We next verify the invariance of the limiting density.
\[
\|\mathcal L_\omega \rho_\omega\|_\infty
\le\|\mathcal L_\omega\mathbf1_\omega\|_\infty\|\rho_\omega\|_\infty
\le C_0\|\rho_\omega\|_\infty,
\qquad
\mathcal L_\omega \rho_{\omega,n}=\rho_{\theta\omega,n+1}.
\]
The first estimate follows from $\mathcal L_\omega\mathbf1_\omega\in C_A$ and~\eqref{eq:cone-bounds}. It permits passage to the uniform limit in the second identity, giving
$\mathcal L_\omega \rho_\omega=\rho_{\theta\omega}$
for every $\omega\in\Omega$.

We verify measurability on $\Omega$ directly.
Use the measurable increasing real lifts $\phi_\omega$
and their iterates $\Phi$ introduced in
Appendix~\ref{app:topological-linearization}.
Fix $n\ge1$ and a real representative $x$ of a circle point.
For each integer $0\le j<d^n$, define
\[
z_j(\omega)
=\Phi(n,\theta^{-n}\omega,\cdot)^{-1}(x+j).
\]
The $n$-step lift $\Phi(n,\theta^{-n}\omega,\cdot)$ is measurable in the base point $\omega$ at each fixed value of the real state variable, and
\[
\{\omega:z_j(\omega)<r\}
=\{\omega:x+j<\Phi(n,\theta^{-n}\omega,r)\},
\qquad r\in\mathbb Q.
\]
Thus $z_j$ is measurable in the base point $\omega$.
For each fixed value $t$ of the real state variable, $D\Phi(n,\theta^{-n}\omega,t)$ is measurable in the base point $\omega$ as a limit of difference quotients. For each fixed $\omega$, it is continuous in the real state variable $t$.
Approximating $z_j$ by rational-valued measurable functions shows that $D\Phi(n,\theta^{-n}\omega,z_j(\omega))$ is measurable.
The indices $0\le j<d^n$ represent exactly all preimages on the circle. Hence the finite sum
\[
\rho_{\omega,n}(x)=\sum_{j=0}^{d^n-1}
\frac{1}{D\Phi(n,\theta^{-n}\omega,z_j(\omega))}
\]
is measurable in the base point $\omega$.
Its limit $\rho_\omega(x)$ is therefore also measurable in the base point $\omega$ for each fixed $x$.

Continuity in the state variable for each fixed $\omega$ and measurability in the base point $\omega$ for each fixed $x$ give joint measurability of $(\omega,x)\mapsto \rho_\omega(x)$.

For uniqueness we first extend the memory-loss estimate~\eqref{eq:cone-memory-loss} to signed functions. Fix $\omega\in\Omega$ and let $\widetilde\rho_\omega\in C^1({\mathbb{S}^1})$ be real-valued with integral zero; this auxiliary function is not a probability density.
If $\widetilde\rho_\omega\ne0$, define the probability density
\[
\widehat\rho_\omega=\mathbf1_\omega
 +\frac{\widetilde\rho_\omega}{2\|\widetilde\rho_\omega\|_{C^1}}.
\]
Both $\widehat\rho_\omega$ and $\mathbf1_\omega$ belong to $C_A$:
the first has integral one, $\widehat\rho_\omega\ge1/2$, and $|D\widehat\rho_\omega|\le1/2\le A\widehat\rho_\omega$, since $A\ge1$.
Using
\[
\widetilde\rho_\omega
=2\|\widetilde\rho_\omega\|_{C^1}
 (\widehat\rho_\omega-\mathbf1_\omega)
\]
in~\eqref{eq:cone-memory-loss} gives
\begin{equation}\label{eq:direct-C1-decay}
\|\mathcal L_\omega^n\widetilde\rho_\omega\|_{C^1}
\le2C_*\|\widetilde\rho_\omega\|_{C^1}
 (1-\tau)^{\lfloor n/N\rfloor}.
\end{equation}
This also holds for $\widetilde\rho_\omega=0$.

To remove the $C^1$ restriction, fix a real-valued zero-mean
$\widetilde\rho_\omega\in L^1({\mathbb{S}^1},m)$ and choose smooth zero-mean
$\widetilde\rho_{\omega,j}$ with $\|\widetilde\rho_\omega-\widetilde\rho_{\omega,j}\|_1\to0$
by periodic convolution.
This approximation is performed on the fixed fiber and requires no measurability of this choice in the base point $\omega$.
The $L^1$ contraction and~\eqref{eq:direct-C1-decay} imply
\[
\begin{aligned}
\|\mathcal L_\omega^n\widetilde\rho_\omega\|_1
&\le\|\mathcal L_\omega^n(\widetilde\rho_\omega-\widetilde\rho_{\omega,j})\|_1
 +\|\mathcal L_\omega^n\widetilde\rho_{\omega,j}\|_1\\
&\le\|\widetilde\rho_\omega-\widetilde\rho_{\omega,j}\|_1
 +2C_*\|\widetilde\rho_{\omega,j}\|_{C^1}
       (1-\tau)^{\lfloor n/N\rfloor}.
\end{aligned}
\]
For each fixed $j$, taking the upper limit as $n\to\infty$
bounds the left-hand side by $\|\widetilde\rho_\omega-\widetilde\rho_{\omega,j}\|_1$.
Letting $j\to\infty$ proves
\begin{equation}\label{eq:direct-L1-decay}
\|\mathcal L_\omega^n\widetilde\rho_\omega\|_1\longrightarrow0
\quad\text{for every }\omega\in\Omega
\text{ and every zero-mean }\widetilde\rho_\omega\in L^1({\mathbb{S}^1},m).
\end{equation}

The invariant-cone and pullback construction above follows \cite[Proposition~5.4]{Stadlbauer}.
We use preservation of $\mathbb P$ to turn the fiberwise decay
in~\eqref{eq:direct-L1-decay} into uniqueness among absolutely continuous
random invariant probability measures.

Let $\bar\mu_\omega$ be such a random probability measure. On a measurable $\mathbb P$-null set containing the exceptional fibers for its probability, absolute-continuity and invariance properties, set $\bar\mu_\omega=m$.
Define its total variation distance from $\mu_\omega$ by
\[
R(\omega)=\|\bar\mu_\omega-\mu_\omega\|_{\mathrm{TV}}
:=\sup_{\substack{v\in C({\mathbb{S}^1})\\\|v\|_\infty\le1}}
\left|\int_{{\mathbb{S}^1}} v\,d\bar\mu_\omega-\int_{{\mathbb{S}^1}} v\,d\mu_\omega\right|.
\]
Here $C({\mathbb{S}^1})$ denotes the space of continuous real-valued functions on ${\mathbb{S}^1}$.
The unit ball of $C({\mathbb{S}^1})$ has a countable dense subset, and integration against either probability measure is continuous in the uniform norm. This countable subset therefore gives the same supremum. Each integral in the resulting countable supremum is measurable, so $R$ is measurable and $0\le R\le2$.

Since $\theta$ preserves $\mathbb P$, for almost every $\omega$
the absolute continuity and invariance relations for the random
measure $\bar\mu_\omega$ hold at $\theta^n\omega$ for every $n\ge0$.
For each such $\omega$, choose a Radon--Nikodym density
$\bar\rho_\omega=d\bar\mu_\omega/dm$ on $\mathbb S^1$.
Then $\bar\rho_\omega-\rho_\omega\in L^1(m)$ has integral zero.
Changing variables on the inverse branches in~\eqref{eq:transfer-operator} shows that the density of a pushed-forward measure is its transfer-operator image.
Thus, for every $n\ge0$,
\[
\bar\mu_{\theta^n\omega}
 =(\mathcal L_\omega^n\bar\rho_\omega)m,
\quad
\mu_{\theta^n\omega}
 =(\mathcal L_\omega^n\rho_\omega)m.
\]
Indeed, random invariance of both $\bar\mu_\omega$ and $\mu_\omega$ and the identity
$(\psi_\omega)_*(f\,m)=(\mathcal L_\omega f)m$, obtained by
change of variables on the inverse branches, give
\[
\bar\mu_{\theta^n\omega}-\mu_{\theta^n\omega}
=\bigl(\mathcal L_\omega^n(\bar\rho_\omega-\rho_\omega)\bigr)m.
\]
For every real $g\in L^1(m)$, the variation measure of $g\,m$ is
$|g|\,m$, so $\|g\,m\|_{\mathrm{TV}}=\int |g|\,dm$.
The formula $\|g\,m\|_{\mathrm{TV}}=\int |g|\,dm$ agrees with the supremum over continuous test functions used in the definition of $R$: continuous functions are dense in
$L^1(|g|m)$, and truncation to $[-1,1]$ gives continuous
approximations to $\operatorname{sgn}(g)$ with absolute value at most one.
Here $\operatorname{sgn}(t)$ equals $1$, $0$ or $-1$ according as $t>0$, $t=0$ or $t<0$.
Consequently,
\[
\begin{aligned}
R(\theta^n\omega)
&=\|\bar\mu_{\theta^n\omega}-\mu_{\theta^n\omega}\|_{\mathrm{TV}}=\int_{{\mathbb{S}^1}}
  \left|\mathcal L_\omega^n(\bar\rho_\omega-\rho_\omega)\right|\,dm=\|\mathcal L_\omega^n(\bar\rho_\omega-\rho_\omega)\|_1
\longrightarrow0.
\end{aligned}
\]
The limit follows from~\eqref{eq:direct-L1-decay}, since
$\bar\rho_\omega-\rho_\omega$ is real, integrable, and has integral zero.
Since $0\le R\circ\theta^n\le2$, dominated convergence and
preservation of $\mathbb P$ give
\[
\int_\Omega R\,d\mathbb P
=\int_\Omega R\circ\theta^n\,d\mathbb P\longrightarrow0.
\]
Hence $R=0$ almost surely, proving uniqueness.
This completes the proof of Lemma 3.
\end{proof}
The density $\bar\rho_\omega$ was chosen only on each fixed fiber; it need not be smooth or belong to the cone.
No measurable selection of these representatives and no ergodicity assumption on the base transformation $\theta$ are required.

\subsection{Higher-order smoothness of the random density}

All inverse branches below are defined on an open arc $I\subset {\mathbb{S}^1}$ of length less than $1/2$.
Let $\Psi_k(-n,\omega,\cdot)$, $1\leq k\leq {d^n}$, denote the local inverse branches of $\Psi(n,\theta^{-n}\omega,\cdot)$, and write $x_{n,k}=\Psi_k(-n,\omega,x)$, $y_{n,k}=\Psi_k(-n,\omega,y)$.
Thus $\Psi(n,\theta^{-n}\omega,x_{n,k})=x$ and similarly for $y$. Set $J_{n,\omega}=D\Psi(n,\theta^{-n}\omega,\cdot)>0$.
Branch labels differ only by permutations on overlaps. Uniform short-arc estimates for sums over all branches therefore extend to ${\mathbb{S}^1}$, with antipodal pairs covered by continuity in the state variable.

By Lemma~\ref{lm4.1}, the pullback densities $\rho_{\omega,n}$ converge to $\rho_\omega$ in the $C^1$ norm in the state variable, uniformly over the base points $\omega\in\Omega$.
We now establish $C^{K-1,\alpha}$ bounds in the state variable for these approximations, uniform over $n$ and the base points $\omega$, and pass them to the limit.
We retain the notation $\rho_{\omega,n}$ for every finite pullback density on $\mathbb{S}^1$ at the base point $\omega$ and write $D^i\rho_{\omega,n}$ for its $i$-th derivative with respect to the state variable. The inverse-branch formula~\eqref{eq:transfer-operator} gives
\[
\rho_{\omega,n}(x)=\sum_{k=1}^{{d^n}}\frac1{J_{n,\omega}(x_{n,k})}.
\]

The following lemma adapts the higher-order density estimates of
Arteaga \cite[Proposition~2.1]{Arteaga1994} to the random cocycle setting.
The proof follows the same differentiation and distortion arguments, with
deterministic iterates replaced by the pullback cocycle iterates
$\Psi(n,\theta^{-n}\omega,\cdot)$ and all estimates uniform in the base
point $\omega$ and the iteration time $n$.
We postpone the lengthy derivative calculations to
Appendix~\ref{app:higher-order-density}.

For a cocycle map $\Psi(n,\omega,\cdot)$, the parenthesized superscript $(r)$ denotes its $r$-th derivative with respect to the state variable.

\begin{lm}\label{lm:higher-order-estimates}
With the notation above, the following assertions hold.
\begin{enumerate}
\renewcommand{\labelenumi}{(\roman{enumi})}
\item
For every integer $1\leq i\leq K-1$, we have
\begin{align}\label{lmeq4-9-1}
D^i\rho_{\omega,n}(x) &=
\sum_{k=1}^{{d^n}}
\left(
\sum_{s=2}^{i}
\sum_{\substack{
    2\leq i_1,\ldots,i_s\leq i\\
    i_1+\cdots+i_s=i+s
}}
C_{i_1,\ldots,i_s}
\frac{
    \displaystyle\prod_{v=1}^{s}
    \bigl(\Psi(n,\theta^{-n}\omega,\cdot)\bigr)^{(i_v)}
    (x_{n,k})
}{
    \bigl(J_{n,\omega}(x_{n,k})\bigr)^{i_1+\cdots+i_s+1}
}\right.
\nonumber\\
&\quad\left.-\frac{
    \bigl(\Psi(n,\theta^{-n}\omega,\cdot)\bigr)^{(i+1)}
    (x_{n,k})
}{
    \bigl(J_{n,\omega}(x_{n,k})\bigr)^{i+2}
}
\right),
\end{align}
where $C_{i_1,\ldots,i_s}=C_{i_1,\ldots,i_s}(i)\in\mathbb{Z}$ depends only on $i$ and the multi-index $(i_1,\ldots,i_s)$.
The sum over $s$ is understood to be zero when $i=1$.

\item
For every integer $1\leq i\leq K$, there exists a constant
$C_1>0$, independent of $\omega$ and $x$, such that
\begin{equation}\label{lmeq4-10-1}
\sup_{\substack{n\geq 1\\1\leq k\leq {d^n}}}
\frac{
    \left|
    \bigl(\Psi(n,\theta^{-n}\omega,\cdot)\bigr)^{(i)}
    (x_{n,k})
    \right|
}{
    \bigl(J_{n,\omega}(x_{n,k})\bigr)^i
}
\leq C_1.
\end{equation}
The derivative ratios in~\eqref{lmeq4-10-1} divide the $i$-th derivative by the $i$-th power of the first derivative of the same cocycle map.

\item
For a fixed integer $K\geq2$ and $n\geq1$, let
\[
G_{n,K,\omega}(x)
=
\sum_{k=1}^{{d^n}}
\frac{D^K\Psi(n,\theta^{-n}\omega,x_{n,k})}
     {J_{n,\omega}(x_{n,k})^{K+1}}.
\]
There is a constant $C_H<\infty$, independent of $n$ and $\omega$, such that
\begin{equation}\label{eq:higher-order-holder}
|G_{n,K,\omega}(x)-G_{n,K,\omega}(y)|
\le C_H d(x,y)^\alpha,
\qquad x,y\in\mathbb S^1.
\end{equation}
In particular, these functions are equicontinuous in the state variable.
\end{enumerate}
\end{lm}

The preceding derivative and H\"older estimates provide the compactness
needed to transfer higher regularity to the limit function.

\begin{lm}\label{lm4.7}
For every base point $\omega\in\Omega$, the density $\rho_\omega$ has $C^{K-1,\alpha}$ regularity in the state variable.
There is a constant $C_7<\infty$, independent of $n$ and $\omega$, such that
\[
\sup_{\omega\in\Omega,\,n\ge0}\|\rho_{\omega,n}\|_{C^{K-1,\alpha}}\le C_7,
\qquad
\sup_{\omega\in\Omega}\|\rho_\omega\|_{C^{K-1,\alpha}}\le C_7.
\]
The convergence $\rho_{\omega,n}\to \rho_\omega$ in the $C^{K-1}$ norm in the state variable is uniform over the base points $\omega\in\Omega$.
Here, for $v\in C^{K-1}(\mathbb S^1)$, we use the norm
\begin{center}
$\displaystyle
\|v\|_{C^{K-1}}:=\sum_{j=0}^{K-1}\|D^jv\|_\infty,
\qquad D^0v:=v.$
\end{center}
\end{lm}

\begin{proof}
We first establish uniform $C^{K-1,\alpha}$ bounds for the
approximating densities. For $n\ge1$, the inverse-branch formula
and Lemma~\ref{lm4.1} give
\[
\sum_{k=1}^{d^n}\frac{1}{J_{n,\omega}(x_{n,k})}
=\rho_{\omega,n}(x)\le C_0.
\]
The same lemma provides uniform bounds for the derivatives of
orders zero and one. For $2\le i\le K-1$, parts~(i) and~(ii) of
Lemma~\ref{lm:higher-order-estimates} imply
\[
|D^i\rho_{\omega,n}(x)|
\le C\sum_{k=1}^{d^n}
       \frac{1}{J_{n,\omega}(x_{n,k})}
\le CC_0,
\]
where $C$ is independent of $n$ and $\omega$.

To estimate the highest-order H\"older seminorm, write
\[
D^{K-1}\rho_{\omega,n}
=R_{\omega,n}-G_{n,K,\omega},
\]
where $R_{\omega,n}$ is the sum of the product terms with
$s\ge2$ in~\eqref{lmeq4-9-1}, with $i=K-1$.
When $K=2$, this sum is empty and $R_{\omega,n}=0$.
For $K\ge3$, each factor in these products involves a derivative
of $\Psi$ of order at most $K-1$. Differentiating with respect
to $x$, and using
$D x_{n,k}=1/J_{n,\omega}(x_{n,k})$, therefore requires derivatives
of $\Psi$ only through order $K$.
Each resulting term is a product of the derivative ratios
bounded in~\eqref{lmeq4-10-1}, multiplied by
$1/J_{n,\omega}(x_{n,k})$.
Since the number of product types depends only on $K$, we obtain
\[
|D R_{\omega,n}(x)|
\le C\sum_{k=1}^{d^n}
       \frac{1}{J_{n,\omega}(x_{n,k})}
\le CC_0.
\]
The branch sums agree on overlapping arcs, so this gives a
uniform Lipschitz bound on the circle.
Combining this bound with part~(iii) of
Lemma~\ref{lm:higher-order-estimates}, and using
$d(x,y)\le d(x,y)^\alpha$, yields
\[
|D^{K-1}\rho_{\omega,n}(x)
 -D^{K-1}\rho_{\omega,n}(y)|
\le (CC_0+C_H)d(x,y)^\alpha.
\]
Together with the derivative bounds above and
$\rho_{\omega,0}=\mathbf1_\omega$, this proves that there is
a constant $C_7$, independent of $n$ and $\omega$, such that
$
\|\rho_{\omega,n}\|_{C^{K-1,\alpha}}\le C_7
~(n\ge0,\ \omega\in\Omega).
$

We next transfer these bounds to the limiting density.
Fix $\omega\in\Omega$.
The lower-order derivatives are uniformly Lipschitz by the
bounds on their next derivatives, while the derivatives of
order $K-1$ have the uniform H\"older bound just proved.
Arzel\`a--Ascoli therefore gives a common subsequence $(n_j)$
and continuous functions $u_0,\ldots,u_{K-1}$ such that
\[
D^i\rho_{\omega,n_j}\longrightarrow u_i
\quad\text{uniformly on }\mathbb S^1,
\qquad 0\le i\le K-1.
\]
Lemma~\ref{lm4.1} identifies $u_0=\rho_\omega$.
On any real coordinate interval, the fundamental theorem
of calculus gives
\[
D^{i-1}\rho_{\omega,n_j}(x)
-D^{i-1}\rho_{\omega,n_j}(x_0)
=\int_{x_0}^{x}D^i\rho_{\omega,n_j}(t)\,dt,
\qquad 1\le i\le K-1.
\]
Passing to the limit yields
\[
u_{i-1}(x)-u_{i-1}(x_0)
=\int_{x_0}^{x}u_i(t)\,dt.
\]
Starting from $u_0=\rho_\omega$, these identities show
successively that $u_i=D^i\rho_\omega$ for every
$1\le i\le K-1$.
Passing to the limit in the highest-order H\"older estimate
then gives $\rho_\omega\in C^{K-1,\alpha}$.
The supremum norms converge along this subsequence, and the
H\"older seminorm is lower semicontinuous under uniform
convergence. Hence
$
\|\rho_\omega\|_{C^{K-1,\alpha}}\le C_7.
$
Since $\omega$ was arbitrary, this bound holds on every fiber.

We now prove that convergence in $C^{K-1}$ is uniform over
$\omega\in\Omega$.
Otherwise, there would exist $\eta>0$, integers $n_j\to\infty$
and base points $\omega_j\in\Omega$ such that, with
$
v_j:=\rho_{\omega_j,n_j}-\rho_{\omega_j},
$
we have $\|v_j\|_{C^{K-1}}\ge\eta$ for every $j$.
The bounds already proved give
$\|v_j\|_{C^{K-1,\alpha}}\le2C_7$.
By Arzel\`a--Ascoli, a subsequence converges in $C^{K-1}$
to a function $v$.
On the other hand, the uniform $C^1$ convergence in
Lemma~\ref{lm4.1} gives
\[
\|v_j\|_{C^1}
\le
\sup_{\omega\in\Omega}
\|\rho_{\omega,n_j}-\rho_\omega\|_{C^1}
\longrightarrow0.
\]
Thus $v=0$, contradicting $\|v_j\|_{C^{K-1}}\ge\eta$.
This proves the stated uniform $C^{K-1}$ convergence.

Finally, we verify joint measurability of the state derivatives.
The function $(\omega,x)\mapsto\rho_\omega(x)$ is jointly
measurable by Lemma~\ref{lm4.1}.
The derivatives established above satisfy, for $1\le i\le K-1$,
\[
D^i\rho_\omega(x)
=\lim_{q\to\infty}
q\bigl(D^{i-1}\rho_\omega(x+q^{-1})
      -D^{i-1}\rho_\omega(x)\bigr),
\]
where $D^0\rho_\omega=\rho_\omega$ and addition is taken
modulo one.
Induction on $i$ shows that each difference quotient,
and hence its pointwise limit, is jointly measurable in
$(\omega,x)$.
The choice of subsequences above does not affect the joint
measurability of $\rho_\omega$ or its state derivatives,
because $\rho_\omega$ is already jointly measurable and its
derivatives are obtained successively as pointwise limits
of jointly measurable difference quotients.
This completes the proof of Lemma 5.
\end{proof}

\begin{proof}[Proof of Proposition~\ref{pro1}]
Lemma~\ref{lm4.1} gives the positive normalized random invariant density $\rho_\omega$, and uniqueness among absolutely continuous random invariant probability measures up to a $\mathbb P$-null set.
Lemma~\ref{lm4.7} gives $\rho_\omega\in C^{K-1,\alpha}$ for every $\omega\in \Omega$, with the uniform bound in~\eqref{eq:uniform-density-regularity}.
This completes the proof of Proposition 2.
\end{proof}

\section{Examples of the two rigidity criteria}\label{sec:examples}
The following examples directly verify the absolute-continuity and
periodic Lyapunov criteria in Theorems~\ref{th2.2}
and~\ref{th:lyapunov-bootstrap}, respectively.

\subsection{An example satisfying the absolute continuity criterion}
\label{subsec:example-explicit-conjugacy}

Let $\Omega=\{0,1\}^{\mathbb Z}$, equipped with the product $\sigma$-algebra and the Bernoulli product probability $\mathbb P$ assigning probability $1/2$ to each symbol.
Write $\omega=(\omega_j)_{j\in\mathbb Z}$, where $\omega_j\in\{0,1\}$.
Let $\theta$ be the left shift, $(\theta\omega)_j=\omega_{j+1}$.
Then $\theta$ is invertible and preserves $\mathbb P$.

Choose $0<t<1/3$. Let $q_0=\mathrm{Id}$, and let $q_1$ be the circle diffeomorphism induced by the increasing degree-one lift
\[
\widehat q_1(x)=x+\frac{t}{2\pi}\sin(2\pi x).
\]
Consider the random expanding circle map
\begin{equation}\label{eq:example-four-maps}
\psi_\omega=q_{\omega_1}^{-1}\circ E_2\circ q_{\omega_0}.
\end{equation}
The pair $(\omega_0,\omega_1)$ determines one of four maps:
\[
\psi_\omega=
\begin{cases}
E_2, & (\omega_0,\omega_1)=(0,0),\\
q_1^{-1}\circ E_2, & (\omega_0,\omega_1)=(0,1),\\
E_2\circ q_1, & (\omega_0,\omega_1)=(1,0),\\
q_1^{-1}\circ E_2\circ q_1,
  & (\omega_0,\omega_1)=(1,1).
\end{cases}
\]

We first verify the fiber-map hypotheses of Theorem~\ref{th2.2}. Put $\widehat q_0(x)=x$ and choose the real lift
\[
\phi_\omega(x)=\widehat q_{\omega_1}^{-1}
   \bigl(2\widehat q_{\omega_0}(x)\bigr).
\]
For each fixed $x$, $\phi_\omega(x)$ is measurable in $\omega$,
and $\phi_\omega(x+1)=\phi_\omega(x)+2$ gives degree two.
Since $1-t\le Dq_i\le1+t$ for $i=0,1$, differentiation gives
\[
D\psi_\omega(x)
=\frac{2Dq_{\omega_0}(x)}{Dq_{\omega_1}(\psi_\omega(x))}
\ge\frac{2(1-t)}{1+t}=: \lambda>1.
\]
The derivative and H\"older bounds are uniform for each
$K\ge2$ and $0<\alpha\le1$, since only four analytic maps occur.

To verify the absolute-continuity criterion, define
$
h_\omega=q_{\omega_0}.
$
Then
\[
h_{\theta\omega}\circ\psi_\omega
=q_{\omega_1}\circ q_{\omega_1}^{-1}
 \circ E_2\circ q_{\omega_0}
=E_2\circ h_\omega.
\]
The bounds $1-t\le Dq_i\le1+t$ make $h_\omega$ and its inverse Lipschitz,
so both preserve null sets and (A1) holds. Moreover,
\[
\nu_\omega:=(h_\omega^{-1})_*m
=Dq_{\omega_0}\,m
=\bigl(1+t\omega_0\cos(2\pi x)\bigr)m.
\]
The density $(\omega,x)\mapsto1+t\omega_0\cos(2\pi x)$ is jointly measurable, so (A2) holds.
Theorem~\ref{th2.2} therefore gives $C^{K,\alpha}$ regularity
of the prescribed conjugacy $h^0_\omega$ for $\mathbb P$-almost every $\omega$.

\subsection{An example satisfying the periodic Lyapunov criterion}
\label{subsec:example-direct-exponent}

Let the base space be $\Omega={\mathbb{S}^1}$ with normalized Lebesgue probability
$\mathbb P$, and let $\theta\omega=\omega+\beta\pmod1$, where $\beta\in\mathbb R\setminus\mathbb Q$.
This transformation is invertible and preserves $\mathbb P$.
Fix $0<c<1/2$ and consider the degree-two system
\begin{equation}\label{eq:example-irrational-affine}
\psi_\omega(x)=2x+\gamma(\omega)\pmod1,
\qquad
\gamma(\omega)=\frac12+c\sin(2\pi\omega).
\end{equation}
These are degree-two maps with $D\psi_\omega=2$ and
$D^j\psi_\omega=0$ for $j\ge2$. Thus the fiber-map hypotheses of Theorem~\ref{th:lyapunov-bootstrap} hold
with $\lambda=2$ for every $K\ge2$ and $0<\alpha\le1$.

For every $\omega$, $x$ and $n\ge1$, we have $D\psi_\omega(x)=2$ and $D\Psi(n,\omega,x)=2^n$.
Hence, for every random periodic orbit $Y$ and each phase $i$,
\[
\frac1n\log D\Psi(n,\omega,Y_i(\omega))=\log2.
\]
Theorem~\ref{th:lyapunov-bootstrap} gives $C^{K,\alpha}$ regularity
of the prescribed conjugacy $h^0_\omega$ for $\mathbb P$-almost every $\omega$.

To exhibit the random periodic orbits explicitly, define
\begin{equation}\label{eq:example-rotation-series}
b_\omega=\sum_{j=0}^{\infty}
                \frac{\gamma(\theta^j\omega)}{2^{j+1}}.
\end{equation}
The series converges uniformly to a bounded measurable function satisfying
$
2b_\omega=\gamma(\omega)+b_{\theta\omega}.
$
For the particular choice of $\gamma$, summing a geometric series gives
\[
b_\omega=\frac12+c\,\operatorname{Im}
             \left(\frac{e^{2\pi\mathrm i\omega}}
             {2-e^{2\pi\mathrm i\beta}}\right).
\]
Here $\operatorname{Im}$ denotes the imaginary part and $\mathrm i^2=-1$.
Consequently $x\mapsto x+b_\omega$ is a real conjugacy from $x\mapsto2x+\gamma(\omega)$ to $x\mapsto2x$.
Since $0<\gamma(\omega)<1$, the lift $x\mapsto2x+\gamma(\omega)$ is already normalized as in the proof of Proposition~\ref{th:continuous-linearization}.
The real conjugacy $x\mapsto x+b_\omega$ has bounded displacement, so uniqueness for the fixed real equation~\eqref{eq:topological-lift-conjugacy} identifies $x\mapsto x+b_\omega$ with the prescribed lift $\widehat h^0_\omega$.
Thus, for $a\in\mathbb S^1$,
\[
h^0_\omega(x)=x+b_\omega\pmod1,
\qquad
\xi_a(\omega)=a-b_\omega\pmod1.
\]

For $\Gamma_a=\Gamma_{\xi_a}$, direct iteration gives
\[
T^n(\omega,\xi_a(\omega))
=(\theta^n\omega,\,2^na-b_{\theta^n\omega}\pmod1),
\qquad
T^n\Gamma_a=\Gamma_{2^na\pmod1}.
\]
Hence $E_2^p a=a$ implies $T^p\Gamma_a=\Gamma_a$, as required by the graph definition.
For example,
\[
T\Gamma_0=\Gamma_0,\qquad
T\Gamma_{1/3}=\Gamma_{2/3},\qquad
T\Gamma_{2/3}=\Gamma_{1/3}.
\]
These give random periodic orbits of periods one and two.
The skew product $T$ itself has no periodic points:
$T^p(\omega,x)=(\omega,x)$ would imply $p\beta\in\mathbb Z$, contradicting irrationality.

\appendix
\makeatletter
\newcommand{\appendixprefix@section}{\appendixname\space}
\renewcommand{\@seccntformat}[1]{%
  \@ifundefined{appendixprefix@#1}{}{\csname appendixprefix@#1\endcsname}%
  \csname the#1\endcsname\quad}
\makeatother

\section{Proof of topological random linearization}\label{app:topological-linearization}
\begin{proof}[Proof of Proposition~\ref{th:continuous-linearization}]
For necessity, suppose that such a random circle homeomorphism $h^0_\omega$ exists and take degrees in~\eqref{eq:topological-conjugacy}. Since the homeomorphisms $h^0_\omega$ and $h^0_{\theta\omega}$ preserve orientation, this gives
\[
\deg\psi_\omega
=\deg(h^0_{\theta\omega}\circ\psi_\omega)
=\deg(E_d\circ h^0_\omega)=d
\qquad (\omega\in\Omega).
\]
For sufficiency, assume $\deg\psi_\omega=d$ for every $\omega\in\Omega$.
Choose increasing real lifts
$\phi_\omega:\mathbb R\to\mathbb R$ of $\psi_\omega$,
measurable in $\omega$ for each fixed real state variable,
and write
\begin{equation}\label{eq:lift-notation}
\phi_\omega(t)=dt+F_\omega(t),
\qquad F_\omega(t+1)=F_\omega(t).
\end{equation}
Such a measurable choice can be made by taking the
representative of $\psi_\omega(0)$ in $[0,1)$ and lifting
the path $s\mapsto\psi_\omega(st\bmod1)$, $0\le s\le1$.
Adding signed increments along dyadic subdivisions to
this representative gives measurable approximations
converging to $\phi_\omega(t)$ for each fixed $t$.
The corresponding real iterates are
\[
\Phi(n,\omega)
=\phi_{\theta^{n-1}\omega}\circ\cdots\circ\phi_\omega,
\qquad n\ge1,
\qquad \Phi(0,\omega)=\mathrm{Id}.
\]
We seek real conjugacies of the form $\mathrm{Id}+g_\omega$
and construct the two directions by the same iteration.
We first normalize the lifts to obtain uniform displacement bounds. Set
\[
\overline F_\omega(x)
=F_\omega(x)-\lfloor F_\omega(0)\rfloor,
\qquad
\overline\phi_\omega(x)
=dx+\overline F_\omega(x).
\]
Subtracting an integer does not change the expanding circle map.
Equal degree and monotonicity give
\begin{equation}\label{eq:normalized-lift-bounds}
\|\overline F_\omega\|_\infty\le d+1,
\qquad D\overline\phi_\omega\ge\lambda>1.
\end{equation}
Indeed, $\overline\phi_\omega(0)=\overline F_\omega(0)\in[0,1)$ and
$\overline\phi_\omega(1)=\overline\phi_\omega(0)+d$.
For $0\le x\le1$, monotonicity bounds $\overline\phi_\omega(x)$ between these two endpoints, so
$|\overline\phi_\omega(x)-dx|\le d+1$; periodicity extends this bound to all real $x$.
Thus no bound on the original constants $F_\omega(0)$ is needed.
Since $d=\int_0^1D\overline\phi_\omega\ge\lambda$, both $\overline\phi_\omega^{-1}$ and $(d\,\mathrm{Id})^{-1}$ are $\varepsilon$-Lipschitz, where $\varepsilon=\lambda^{-1}$.

Let $\mathcal B$ be the space of continuous random functions $g_\omega:\mathbb R\to\mathbb R$ that are measurable in the base point $\omega$ for each fixed $x$, continuous and $1$-periodic in the state variable $x$ for each fixed $\omega$, and satisfy
\[
\|g\|_{\mathcal B}
:=\sup_{\omega\in\Omega,\,x\in\mathbb R}|g_\omega(x)|
<\infty.
\]
This space is complete in the uniform norm.
We write $\mathcal H=\mathrm{Id}+\mathcal B$ for the corresponding degree-one real maps $H_\omega$, continuous in the state variable, with displacement uniformly bounded over $(\omega,x)\in \Omega\times\mathbb R$.

To treat the conjugacy and its inverse together, consider
\[
\varphi_\omega(x)=dx+\widetilde F_\omega(x),
\qquad
\varphi_\omega^*(x)=dx+\widetilde F_\omega^*(x),
\]
where each auxiliary displacement is chosen to be either $\overline F_\omega$ or zero.

We seek a conjugacy from $\varphi_\omega$ to $\varphi_\omega^*$ by solving
\begin{equation}\label{eq:iteration-general-conjugacy}
(\mathrm{Id}+g_{\theta\omega})\circ\varphi_\omega
=\varphi_\omega^*\circ(\mathrm{Id}+g_\omega)
\end{equation}
in $\mathcal H$.
Since $\varphi_\omega^*$ is an increasing homeomorphism of $\mathbb R$, we can define the iteration, starting from $g^{(0)}=0$, by
\begin{equation}\label{eq:iteration-explicit}
\mathrm{Id}+g_\omega^{(n+1)}
=(\varphi_\omega^*)^{-1}\circ
(\mathrm{Id}+g_{\theta\omega}^{(n)})\circ\varphi_\omega.
\end{equation}

For every fixed value $y$ of the real state variable, the inverse $(\varphi_\omega^*)^{-1}(y)$ is measurable in the base point $\omega$, because, for each rational $r$,
\[
(\varphi_\omega^*)^{-1}(y)<r
\quad\Longleftrightarrow\quad
y<\varphi_\omega^*(r).
\]
For each fixed base point $\omega$, the inverse is also continuous in the state variable $y$.
Measurability in the base point for each fixed value of the state variable, together with continuity in that variable, permits evaluation at values of the state variable that depend measurably on the base point, by approximation with countably valued measurable functions.
Moreover, the equal-degree translation rule gives
\[
\varphi_\omega(x+1)=\varphi_\omega(x)+d,
\qquad
(\varphi_\omega^*)^{-1}(y+d)
=(\varphi_\omega^*)^{-1}(y)+1.
\]
Since $d$ is an integer, the $1$-periodicity of
$g_{\theta\omega}^{(n)}$ therefore implies the
$1$-periodicity of $g_\omega^{(n+1)}$.

For $g^{(n)}\in\mathcal B$, the iteration thus gives a displacement that is continuous and $1$-periodic in the state variable $x$ for each fixed $\omega$, and measurable in the base point $\omega$ for each fixed $x$.
Since $(\varphi_\omega^*)^{-1}$ is $\varepsilon$-Lipschitz, comparison of
$
(\varphi_\omega^*)^{-1}
\bigl(
\varphi_\omega(x)
+g_{\theta\omega}^{(n)}(\varphi_\omega(x))
\bigr)
$
with
$
(\varphi_\omega^*)^{-1}(\varphi_\omega^*(x))=x
$
gives
\[
\|g^{(n+1)}\|_{\mathcal B}
\le
\varepsilon\bigl(
\|g^{(n)}\|_{\mathcal B}+2(d+1)
\bigr).
\]
Thus every iterate belongs to $\mathcal B$.

The same inverse Lipschitz estimate applied to \eqref{eq:iteration-explicit} gives, for $n\ge1$,
\[
\|g^{(n+1)}-g^{(n)}\|_{\mathcal B}
\le
\varepsilon
\|g^{(n)}-g^{(n-1)}\|_{\mathcal B}.
\]
Also, $\|g^{(1)}\|_{\mathcal B}\le 2(d+1)/\lambda$.
Consequently,
\begin{equation}\label{eq:iteration-cauchy}
\|g^{(n+1)}-g^{(n)}\|_{\mathcal B}
\le
\frac{2(d+1)}{\lambda^{n+1}},
\qquad n\ge0.
\end{equation}
As $\lambda>1$, the sequence converges in $\mathcal B$ to a function $g$ with $\|g\|_{\mathcal B}\le 2(d+1)/(\lambda-1)$.
Passing to the limit in \eqref{eq:iteration-explicit} proves \eqref{eq:iteration-general-conjugacy}.

For these fixed real lifts $\varphi_\omega$ and $\varphi_\omega^*$, the solution of the exact real equation~\eqref{eq:iteration-general-conjugacy} is unique in $\mathcal H$.
Indeed, if $\mathrm{Id}+q_1$ and $\mathrm{Id}+q_2$ are two such solutions, then
\[
\mathrm{Id}+q_{i,\omega}
=(\varphi_\omega^*)^{-1}\circ
(\mathrm{Id}+q_{i,\theta\omega})\circ\varphi_\omega,
\qquad i=1,2.
\]
The inverse Lipschitz estimate yields
\[
\|q_1-q_2\|_{\mathcal B}
\le
\varepsilon\|q_1-q_2\|_{\mathcal B}.
\]
The finiteness of the norm, together with $\varepsilon<1$, gives $q_1=q_2$.
This is uniqueness for a fixed real equation in $\mathcal H$; circle conjugacies satisfying only equality modulo integers need not be unique.

Take
$(\widetilde F,\widetilde F^*)=(\overline F,0)$
and $(0,\overline F)$ in turn.
The resulting lifts
$\widehat h_1,\widehat h_2\in\mathcal H$ satisfy
\begin{equation}\label{eq:iteration-two-directions}
\begin{aligned}
\widehat h_{1,\theta\omega}\circ\overline\phi_\omega
&=(d\,\mathrm{Id})\circ\widehat h_{1,\omega},\quad
\widehat h_{2,\theta\omega}\circ(d\,\mathrm{Id})
&=\overline\phi_\omega\circ\widehat h_{2,\omega}.
\end{aligned}
\end{equation}
To prove that the two constructed maps are mutually inverse,
note first that $\mathcal H$ is closed under fiberwise composition.
Indeed, if $H_i=\mathrm{Id}+q_i$, then
\[
H_1\circ H_2-\mathrm{Id}
=q_2+q_1\circ(\mathrm{Id}+q_2)\in\mathcal B.
\]
By \eqref{eq:iteration-two-directions},
\[
\begin{aligned}
(\widehat h_{2,\theta\omega}
 \circ\widehat h_{1,\theta\omega})
 \circ\overline\phi_\omega
&=
\overline\phi_\omega\circ
(\widehat h_{2,\omega}\circ\widehat h_{1,\omega}),\\
(\widehat h_{1,\theta\omega}
 \circ\widehat h_{2,\theta\omega})
 \circ(d\,\mathrm{Id})
&=
(d\,\mathrm{Id})\circ
(\widehat h_{1,\omega}\circ\widehat h_{2,\omega}).
\end{aligned}
\]
Hence $\widehat h_2\circ\widehat h_1$ and
$\widehat h_1\circ\widehat h_2$ solve
\eqref{eq:iteration-general-conjugacy} for the pairs
$(\overline F,\overline F)$ and $(0,0)$, respectively.
The map $\mathrm{Id}$ solves each of these two fixed real equations.
Uniqueness in $\mathcal H$ therefore gives
\[
\widehat h_{2,\omega}\circ\widehat h_{1,\omega}
=
\widehat h_{1,\omega}\circ\widehat h_{2,\omega}
=
\mathrm{Id}.
\]
Thus these continuous lifts are mutually inverse
homeomorphisms.
Their degree-one translation rule makes them
orientation preserving.

Projecting $\widehat h_{1,\omega}$ to ${\mathbb{S}^1}$
defines $h^0_\omega$ and gives
\eqref{eq:topological-conjugacy}.
The inverse $(h^0_\omega)^{-1}$ is the projection of $\widehat h_{2,\omega}$.
Measurability of $h^0_\omega$ and $(h^0_\omega)^{-1}$ in the base point, at each fixed value of the state variable, follows from the construction.

The lift
$\widehat h^0_\omega:=\widehat h_{1,\omega}$
is therefore the unique element of $\mathcal H$
satisfying
\begin{equation}\label{eq:topological-lift-conjugacy}
\widehat h^0_{\theta\omega}\circ\overline\phi_\omega
=d\widehat h^0_\omega.
\end{equation}
This fixes the random circle homeomorphism $h^0_\omega$ used throughout the paper without
imposing any absolute continuity condition.
This completes the proof of Proposition 1.
\end{proof}

\section{Density of random periodic points}\label{app:periodic-graph-density}

Let $\mathcal M(\Omega,\mathbb{S}^1)$ be the space of measurable sections, equipped with the supremum distance $d_\infty$ in~\eqref{dis-section}.

\begin{pro}\label{pro:periodic-graph-density}
Let $\theta:\Omega\to\Omega$ be invertible and $\mathbb P$-preserving, with $\theta$ and $\theta^{-1}$ measurable.
Suppose that $\psi_\omega:\mathbb{S}^1\to\mathbb{S}^1$ is an orientation-preserving random $C^1$ expanding circle map of degree $d\ge2$ for every $\omega\in\Omega$ and satisfies \eqref{eq:standing-expansion}.
For every $g\in\mathcal M(\Omega,\mathbb{S}^1)$ and every integer $p\ge1$, there is a measurable section $\xi_p$ such that
\begin{equation}\label{eq:periodic-density-return}
T^p\Gamma_{\xi_p}=\Gamma_{\xi_p},
\qquad
\Psi(p,\omega,\xi_p(\omega))=\xi_p(\theta^p\omega)
\quad (\omega\in\Omega),
\end{equation}
and
\begin{equation}\label{eq:periodic-graph-density-bound}
d_\infty(g,\xi_p)\le\frac{1}{2(\lambda^p-1)}.
\end{equation}
Thus the random periodic points are dense in $\mathcal M(\Omega,\mathbb{S}^1)$ for $d_\infty$, and their values are dense in every fiber.
The prescribed $p$ is a return period, not necessarily the least period.
\end{pro}

\begin{proof}
Fix $p\ge1$ and $g\in\mathcal M(\Omega,\mathbb{S}^1)$.

\noindent{\emph{Step 1. Constructing measurable inverses of the real lifts and establishing contraction.}}
\par\nobreak\smallskip\nobreak
Choose the increasing real lift $\phi_\omega$ with
$\phi_\omega(0)\in[0,1)$ as in
Appendix~\ref{app:topological-linearization}.
It is continuous in the real state variable and measurable in $\omega$ at each fixed value of that variable, hence jointly measurable.
It satisfies
\[
\phi_\omega(t+1)=\phi_\omega(t)+d,\qquad
D\phi_\omega(t)\ge\lambda.
\]
Put
\[
A_\omega:=\Phi(p,\omega)
=\phi_{\theta^{p-1}\omega}\circ\cdots\circ\phi_\omega.
\]
Then $A_\omega(t+1)=A_\omega(t)+d^p$ and $DA_\omega\ge\lambda^p$.
In particular, $A_\omega$ is an increasing bijection of $\mathbb R$ onto itself and
\begin{equation}\label{eq:periodic-density-inverse}
|A_\omega^{-1}(s)-A_\omega^{-1}(t)|
\le\lambda^{-p}|s-t|.
\end{equation}
The inverse $A_\omega^{-1}(t)$ is jointly measurable in $(\omega,t)$: for every rational $r$,
\[
A_\omega^{-1}(t)<r
\quad\Longleftrightarrow\quad
t<A_\omega(r),
\]
and the latter condition is measurable in $(\omega,t)$.
Only the real lift is inverted here; the expanding circle map remains a non-invertible covering.

\noindent{\emph{Step 2. Choosing an integer correction to control the section's return error.}}
\par\nobreak\smallskip\nobreak
Let $\widetilde g(\omega)\in[0,1)$ be the measurable representative of $g(\omega)$ and define
\[
\ell(\omega)=\left\lfloor
A_\omega(\widetilde g(\omega))
-\widetilde g(\theta^p\omega)+\frac12
\right\rfloor.
\]
This integer-valued function is measurable.
The residual
\[
r(\omega):=A_\omega(\widetilde g(\omega))
-\widetilde g(\theta^p\omega)-\ell(\omega)
\]
satisfies $|r(\omega)|\le1/2$ for every $\omega\in\Omega$.

\noindent{\emph{Step 3. Solving the return equation by contraction and estimating the correction.}}
\par\nobreak\smallskip\nobreak
Let $\mathcal B_b$ be the Banach space of bounded measurable real functions on $\Omega$, with the actual supremum norm over $\Omega$.
Define
\[
(\mathcal Qb)(\omega)
=A_\omega^{-1}\bigl(\widetilde g(\theta^p\omega)
+\ell(\omega)+b(\theta^p\omega)\bigr)-\widetilde g(\omega).
\]
Joint measurability of the inverse shows that $\mathcal Qb$ is measurable.
By~\eqref{eq:periodic-density-inverse} and the residual identity,
\[
\|\mathcal Qb-\mathcal Qc\|_\infty
\le\lambda^{-p}\|b-c\|_\infty,
\qquad
\|\mathcal Qb\|_\infty
\le\lambda^{-p}\left(\|b\|_\infty+\frac12\right).
\]
Thus $\mathcal Q$ maps $\mathcal B_b$ into itself and is a strict contraction.
Its fixed point $b_*$ is measurable and satisfies
\[
\|b_*\|_\infty
\le\lambda^{-p}\left(\|b_*\|_\infty+\frac12\right),
\qquad
\|b_*\|_\infty\le\frac{1}{2(\lambda^p-1)}.
\]
The fixed-point equation is
\[
A_\omega(\widetilde g(\omega)+b_*(\omega))
=\widetilde g(\theta^p\omega)+b_*(\theta^p\omega)+\ell(\omega).
\]

\noindent{\emph{Step 4. Constructing random periodic graphs and proving density in the section space.}}
\par\nobreak\smallskip\nobreak
Set
\[
\xi_p(\omega)=[\widetilde g(\omega)+b_*(\omega)]\in\mathbb{S}^1{,\qquad\omega\in\Omega}.
\]
{Here $[t]$ denotes the class of $t$ modulo one.}
Projection modulo one removes $\ell(\omega)$ and yields the return relation in~\eqref{eq:periodic-density-return}.
{Since $\theta^p$ is bijective, the section return relation in~\eqref{eq:periodic-density-return}
is equivalent to $T^p\Gamma_{\xi_p}=\Gamma_{\xi_p}$.}
{Also, $d(\xi_p(\omega),g(\omega))\le |b_*(\omega)|$,
which proves~\eqref{eq:periodic-graph-density-bound}.}
The measurable choice in $[0,1)$ and the final projection handle the cut point of the circle without requiring a continuous choice of representatives.

For clarity, the associated representing maps are
\[
Y_i(\omega)=\Psi(i,\theta^{-i}\omega,\xi_p(\theta^{-i}\omega)),
\qquad i\in\mathbb N_0.
\]
{Here $\mathbb N_0=\{0,1,2,\ldots\}$.}
Using the cocycle identity and the return relation at $\theta^{-(i+p)}\omega$ gives
\[
\begin{aligned}
Y_{i+p}(\omega)
&=\Psi\bigl(i,\theta^{-i}\omega,
       \Psi(p,\theta^{-(i+p)}\omega,\xi_p(\theta^{-(i+p)}\omega))\bigr)\\
&=\Psi(i,\theta^{-i}\omega,\xi_p(\theta^{-i}\omega))
=Y_i(\omega).
\end{aligned}
\]
{The cocycle identity and the definition of $Y_i$ also give
$\psi_\omega(Y_i(\omega))=Y_{i+1}(\theta\omega)$.}
Equivalently, $\Gamma_{Y_i}=T^i\Gamma_{\xi_p}$ and $\Gamma_{Y_{i+p}}=\Gamma_{Y_i}$.
{Thus $\xi_p=Y_0$ is a random periodic point in Definition~\ref{def:random-periodic-graphs}.}

Given $\epsilon>0$, every sufficiently large $p$ makes the bound in~\eqref{eq:periodic-graph-density-bound} smaller than $\epsilon$, proving uniform density in the section space.
Taking $g$ to be an arbitrary constant section proves density of the values on every fixed fiber.
This completes the proof of Proposition 3.
\end{proof}

The proof uses neither topological conjugacy nor ergodicity. The countable test collection $\{\xi_a:a\in\mathcal P_d\}$ used in Theorem~\ref{th:lyapunov-bootstrap} has a different density statement.
It is fiberwise dense, but need not be dense in the section space.
For example, take the constant system $\psi_\omega=E_2$ and $h^0_\omega=\mathrm{Id}$ over a base with two disjoint sets $B$ and $\Omega\setminus B$ of positive probability.
All induced test sections are constant.
The measurable section equal to $0$ on $B$ and $1/2$ on its complement has {supremum distance} at least $1/4$ from every constant section, by the triangle inequality on $\mathbb{S}^1$.
{The failure of this countable test collection to be dense in the section space does not contradict} Proposition~\ref{pro:periodic-graph-density}, which allows all random periodic points, including sections whose values in the linear coordinates depend on $\omega$.

\section{Uniform estimates of higher-order derivatives for random density }\label{app:higher-order-density}

We retain the notation of Section~\ref{sec:densities}.
The following proof adapts the three estimates in
\cite[Proposition~2.1]{Arteaga1994} to the pullback cocycle,
with constants independent of the base point and the pullback time.

\begin{proof}[Proof of Lemma~\ref{lm:higher-order-estimates}]
{We prove parts (i), (ii) and (iii) in three steps.}

\noindent{\textit{Step 1. Deriving the expansion of the density derivatives.}}

We prove \eqref{lmeq4-9-1} by induction on $i$.
{The case $i=1$ follows from the inverse-branch differentiation formula~\eqref{eq:transfer-derivative},
applied to $\Psi(n,\theta^{-n}\omega,\cdot)$ with constant density $1$.}

Suppose that \eqref{lmeq4-9-1} holds for some $1\leq i<K-1$.
{For every multi-index in the product-sum term of~\eqref{lmeq4-9-1}, where $s\ge2$,
the product rule and the chain rule give}
\begin{align*}
&\frac{d}{dx}
\left(
\frac{
    \displaystyle\prod_{v=1}^{s}
    \bigl(\Psi(n,\theta^{-n}\omega,\cdot)\bigr)^{(i_v)}
    (x_{n,k})
}{
    \bigl(J_{n,\omega}(x_{n,k})\bigr)^{i_1+\cdots+i_s+1}
}
\right) \\
&=
\sum_{v=1}^{s}
\frac{
    \bigl(\Psi(n,\theta^{-n}\omega,\cdot)\bigr)^{(i_v+1)}
    (x_{n,k})
    \displaystyle\prod_{\substack{1\leq b\leq s\\ b\neq v}}
    \bigl(\Psi(n,\theta^{-n}\omega,\cdot)\bigr)^{(i_b)}
    (x_{n,k})
}{
    \bigl(J_{n,\omega}(x_{n,k})\bigr)^{i_1+\cdots+i_s+2}
} \\
&\quad
-
(i_1+\cdots+i_s+1)
\frac{
    \bigl(\Psi(n,\theta^{-n}\omega,\cdot)\bigr)^{(2)}
    (x_{n,k})
    \displaystyle\prod_{b=1}^{s}
    \bigl(\Psi(n,\theta^{-n}\omega,\cdot)\bigr)^{(i_b)}
    (x_{n,k})
}{
    \bigl(J_{n,\omega}(x_{n,k})\bigr)^{i_1+\cdots+i_s+3}
}.
\end{align*}

In each term obtained by differentiating the numerator, the number of factors remains $s$, while the sum of their derivative orders increases from $i+s$ to
$(i+1)+s$. In the term obtained by differentiating the denominator, one additional factor of derivative order $2$ appears.
{Thus the number of factors becomes $s+1$, and the sum of their derivative orders becomes $i+s+2=(i+1)+(s+1)$.}
Consequently, all these terms have the form required by the first sum in \eqref{lmeq4-9-1} with $i$ replaced by $i+1$.

{The last term in~\eqref{lmeq4-9-1} satisfies}
\begin{align*}
&\frac{\mathrm{d}}{\mathrm{d}x}
\left(
-
\frac{
    \bigl(\Psi(n,\theta^{-n}\omega,\cdot)\bigr)^{(i+1)}
    (x_{n,k})
}{
    \bigl(J_{n,\omega}(x_{n,k})\bigr)^{i+2}
}
\right) \\
&=
-
\frac{
    \bigl(\Psi(n,\theta^{-n}\omega,\cdot)\bigr)^{(i+2)}
    (x_{n,k})
}{
    \bigl(J_{n,\omega}(x_{n,k})\bigr)^{i+3}
} \\
&\quad
+
(i+2)
\frac{
    \bigl(\Psi(n,\theta^{-n}\omega,\cdot)\bigr)^{(2)}
    (x_{n,k})
    \bigl(\Psi(n,\theta^{-n}\omega,\cdot)\bigr)^{(i+1)}
    (x_{n,k})
}{
    \bigl(J_{n,\omega}(x_{n,k})\bigr)^{i+4}
}.
\end{align*}
The first term is precisely the last term in \eqref{lmeq4-9-1} at order $i+1$. The second term belongs to its first sum: it has two factors, of derivative orders $2$ and $i+1$, whose sum is $(i+1)+2$.

After collecting terms, all resulting coefficients are integers and depend only on the derivative order and the corresponding multi-index.
Therefore \eqref{lmeq4-9-1} holds at order $i+1$. This completes the induction.
{This proves part (i).}

\medskip\noindent{\textit{Step 2. Establishing uniform bounds for the derivative ratios.}}

For $i=1$, we have
\[
\frac{
    \left|
    \bigl(\Psi(n,\theta^{-n}\omega,\cdot)\bigr)^{(1)}
    (x_{n,k})
    \right|
}{
    J_{n,\omega}(x_{n,k})
}
=1.
\]

Suppose that the estimates have been established for all derivative orders less than $K$.
We describe the induction step at order $K$; the same argument
applies to every intermediate order.
By the induction hypothesis, there exists $C_2\geq 1$ such that
\begin{equation*}
\frac{
    \left|
    \bigl(\Psi(j,\theta^{-n}\omega,\cdot)\bigr)^{(l)}
    (x_{n,k})
    \right|
}{
    \left(
    \bigl(\Psi(j,\theta^{-n}\omega,\cdot)\bigr)^{(1)}
    (x_{n,k})
    \right)^l
}
\leq C_2,
\qquad
1\leq l\leq K-1,\quad 0\leq j\leq n.
\end{equation*}
For $j=0$, this follows from $\Psi(0,\theta^{-n}\omega,\cdot)=\mathrm{Id}$.
For $j\geq1$, it follows from the uniformity of the induction hypothesis in $\omega$ and $x$.

For $0\leq j\leq n-1$, the cocycle identity gives
\[
\Psi(j+1,\theta^{-n}\omega,\cdot)
=
\psi_{\theta^{-(n-j)}\omega}
\circ
\Psi(j,\theta^{-n}\omega,\cdot).
\]
To make the telescoping step explicit, fix $z=x_{n,k}$ and abbreviate
\[
F_j=\Psi(j,\theta^{-n}\omega,\cdot),\qquad
u_j=F_j(z),\qquad a_j=DF_j(z),\qquad
f_j=\psi_{\theta^{j-n}\omega}.
\]
For a multi-index $(s_1,\ldots,s_{K-1})$, put
\[
t=\sum_{\ell=1}^{K-1}s_\ell,\qquad
A_{s_1,\ldots,s_{K-1}}
=\frac{K!}{\prod_{\ell=1}^{K-1}s_\ell!(\ell!)^{s_\ell}}.
\]
The higher-order chain rule for $F_{j+1}=f_j\circ F_j$, together with
$a_{j+1}=Df_j(u_j)a_j$, gives
\begin{align}
&
\frac{
    \bigl(\Psi(n,\theta^{-n}\omega,\cdot)\bigr)^{(K)}
    (x_{n,k})
}{
    \bigl(J_{n,\omega}(x_{n,k})\bigr)^K
}
\notag\\*
&=
\sum_{j=0}^{n-1}
\left(
\frac{
    \bigl(\Psi(j,\theta^{-n}\omega,\cdot)\bigr)^{(1)}
    (x_{n,k})
}{
    J_{n,\omega}(x_{n,k})
}
\right)^{K-1}
\notag\\*
&\quad\times
\Bigg(
\frac{
    \bigl(\psi_{\theta^{-(n-j)}\omega}\bigr)^{(K)}
    \bigl(\Psi(j,\theta^{-n}\omega,x_{n,k})\bigr)
}{
    \bigl(\psi_{\theta^{-(n-j)}\omega}\bigr)^{(1)}
    \bigl(\Psi(j,\theta^{-n}\omega,x_{n,k})\bigr)
}
\notag\\
&\qquad+
\sum_{\substack{
    s_1,\ldots,s_{K-1}\geq0\\
    \sum_{l=1}^{K-1}l s_l=K\\
    2\leq\sum_{l=1}^{K-1}s_l\leq K-1
}}
A_{s_1,\ldots,s_{K-1}}
\frac{
    \bigl(\psi_{\theta^{-(n-j)}\omega}\bigr)^{(t)}
    \bigl(\Psi(j,\theta^{-n}\omega,x_{n,k})\bigr)
}{
    \bigl(\psi_{\theta^{-(n-j)}\omega}\bigr)^{(1)}
    \bigl(\Psi(j,\theta^{-n}\omega,x_{n,k})\bigr)
}
\notag\\*
&\qquad\qquad\times
\prod_{l=1}^{K-1}
\left(
\frac{
    \bigl(\Psi(j,\theta^{-n}\omega,\cdot)\bigr)^{(l)}
    (x_{n,k})
}{
    \left(
    \bigl(\Psi(j,\theta^{-n}\omega,\cdot)\bigr)^{(1)}
    (x_{n,k})
    \right)^l
}
\right)^{s_l}
\Bigg).
\label{lmeq4-10-3}
\end{align}
The coefficients $A_{s_1,\ldots,s_{K-1}}$ are positive integers.
The inner sum is understood to be zero when $K=2$.

By the bounds on the higher derivatives of the random {expanding} circle map $\psi_\omega$ with respect to the state variable, uniform over the base points $\omega$,

for every $2\leq l\leq K$,
\[
\frac{
    \left|
    \bigl(\psi_{\theta^{-(n-j)}\omega}\bigr)^{(l)}
    \bigl(\Psi(j,\theta^{-n}\omega,x_{n,k})\bigr)
    \right|
}{
    \bigl(\psi_{\theta^{-(n-j)}\omega}\bigr)^{(1)}
    \bigl(\Psi(j,\theta^{-n}\omega,x_{n,k})\bigr)
}
\leq
U\varepsilon.
\]
Moreover, the induction hypothesis gives
\[
\prod_{l=1}^{K-1}
\left|
\frac{
    \bigl(\Psi(j,\theta^{-n}\omega,\cdot)\bigr)^{(l)}
    (x_{n,k})
}{
    \left(
    \bigl(\Psi(j,\theta^{-n}\omega,\cdot)\bigr)^{(1)}
    (x_{n,k})
    \right)^l
}
\right|^{s_l}
\leq C_2^{s_1+\cdots+s_{K-1}}.
\]
The set of admissible multi-indices in \eqref{lmeq4-10-3} is finite. Hence
\[
C_3
:=
\sum_{\substack{
    s_1,\ldots,s_{K-1}\geq0\\
    \sum_{l=1}^{K-1}l s_l=K\\
    2\leq\sum_{l=1}^{K-1}s_l\leq K-1
}}
A_{s_1,\ldots,s_{K-1}}
C_2^{s_1+\cdots+s_{K-1}}
<\infty,
\]
and $C_3$ is independent of $n$, $\omega$, $j$, $x$,
and $k$.

Furthermore, the chain rule and the expanding
property imply
\[
0<
\frac{
    \bigl(\Psi(j,\theta^{-n}\omega,\cdot)\bigr)^{(1)}
    (x_{n,k})
}{
    J_{n,\omega}(x_{n,k})
}
\leq
\frac{1}{\varepsilon^{-(n-j)}},
\qquad 0\leq j\leq n-1.
\]
Combining these estimates with
\eqref{lmeq4-10-3}, we obtain
\begin{align*}
&
\frac{
    \left|
    \bigl(\Psi(n,\theta^{-n}\omega,\cdot)\bigr)^{(K)}
    (x_{n,k})
    \right|
}{
    \bigl(J_{n,\omega}(x_{n,k})\bigr)^K
}
\leq
U\varepsilon(1+C_3)
\sum_{j=0}^{n-1}
\varepsilon^{(K-1)(n-j)}
\\
&~~~~~~~~~~~~~~~~~~~~~~~~~~~~~~~~~=
\frac{
    U(1+C_3)\varepsilon^K
}{
    1-\varepsilon^{K-1}
}
\left(1-\varepsilon^{(K-1)n}\right)
\leq
\frac{
    U(1+C_3)\varepsilon^K
}{
    1-\varepsilon^{K-1}
}
<\infty.
\end{align*}
This bound is independent of $n$, $\omega$, $x$, and $k$, and therefore completes the induction.
Taking $C_1$ to dominate the finitely many bounds obtained for $1\leq i\leq K$ establishes part~(ii).

\medskip\noindent{\textit{Step 3. Establishing a uniform H\"older bound for $G_{n,K,\omega}$.}}

{Let $C_2,C_3$ be the constants obtained in the proof of part~(ii).} Take $C_1$ to be a common constant in~\eqref{lmeq4-10-1} for all derivative orders $1,\ldots,K$, increasing it if necessary so that $C_1\geq C_2\geq1$.

{We separate the contribution involving $D^K\psi$ from the lower-order remainder.
We will prove a H\"older bound for the former and a Lipschitz bound for the latter.}

Fix $n,\omega$ and an open arc $I$ of length less than $1/2$. For $x,y\in I$, pair
\[
x_{n,k}=\Psi_k(-n,\omega,x),
\qquad
y_{n,k}=\Psi_k(-n,\omega,y)
\]
through the same inverse branch.
Define the positive factor
\[
A^j_{n,k}(x)
:=
\frac{
\bigl(D\Psi(j,\theta^{-n}\omega,x_{n,k})\bigr)^{K-1}
}{
D\psi_{\theta^{j-n}\omega}\bigl(
          \Psi(j,\theta^{-n}\omega,x_{n,k})\bigr)
J_{n,\omega}(x_{n,k})^K
},
\]
and set
\[
g^j_{x_{n,k}}
:=
D^K\psi_{\theta^{j-n}\omega}\bigl(
            \Psi(j,\theta^{-n}\omega,x_{n,k})\bigr)
A^j_{n,k}(x),
\qquad
g_{n,K,\omega}(x)
:=
\sum_{k=1}^{{d^n}}\sum_{j=0}^{n-1}g^j_{x_{n,k}}.
\]

{Both parts of the decomposition use bounds on $A^j_{n,k}$ and its logarithmic derivative.}

Differentiating the inverse-branch identity
$
\Psi(n,\theta^{-n}\omega,\Psi_k(-n,\omega,x))=x
$
and using expansion gives
\[
\frac{\mathrm{d}x_{n,k}}{\mathrm{d}x}
=
D\Psi_k(-n,\omega,x)
=
\frac1{J_{n,\omega}(x_{n,k})},
\]
and
\[
0<
\frac{D\Psi(j,\theta^{-n}\omega,x_{n,k})}
     {J_{n,\omega}(x_{n,k})}
\leq\varepsilon^{n-j}.
\]
Since $D\psi_\eta\geq\varepsilon^{-1}$ for every $\eta\in{\Omega}$, it follows that
\[
0<A^j_{n,k}(x)
\leq
\frac{\varepsilon\,\varepsilon^{(K-1)(n-j)}}
     {J_{n,\omega}(x_{n,k})}.
\]

The logarithmic derivative of $A^j_{n,k}$ with respect to the state variable $x$ is
\[
\begin{aligned}
\frac{DA^j_{n,k}(x)}{A^j_{n,k}(x)}
={}&
(K-1)
\frac{
D^2\Psi(j,\theta^{-n}\omega,x_{n,k})
}{
D\Psi(j,\theta^{-n}\omega,x_{n,k})
J_{n,\omega}(x_{n,k})
}
\\
&-
\frac{
D^2\psi_{\theta^{j-n}\omega}\bigl(
            \Psi(j,\theta^{-n}\omega,x_{n,k})\bigr)
}{
D\psi_{\theta^{j-n}\omega}\bigl(
          \Psi(j,\theta^{-n}\omega,x_{n,k})\bigr)
}
\frac{D\Psi(j,\theta^{-n}\omega,x_{n,k})}
     {J_{n,\omega}(x_{n,k})}
\\
&-
K\frac{
D^2\Psi(n,\theta^{-n}\omega,x_{n,k})
}{
J_{n,\omega}(x_{n,k})^2
}.
\end{aligned}
\]
Estimate~\eqref{lmeq4-10-1}, the bound $|D^2\psi_\eta|\leq U$ for every $\eta\in{\Omega}$, and expansion imply
\[
\left|
\frac{DA^j_{n,k}(x)}{A^j_{n,k}(x)}
\right|
\leq
\bigl((K-1)C_1+U\varepsilon\bigr)\varepsilon^{n-j}
+KC_1
\leq C_4,
\]
where
$
C_4:=(2K-1)C_1+U\varepsilon.
$

We first estimate the lower-order remainder $G_{n,K,\omega}-g_{n,K,\omega}$.
For $K=2$, this remainder is zero, and $C_3=0$.
Suppose $K\geq3$.
By \eqref{lmeq4-10-3}, each remainder term consists of the same positive factor $A^j_{n,k}$ multiplied by a fixed coefficient and a factor of the form
\[
D^t\psi_{\theta^{j-n}\omega}\bigl(
            \Psi(j,\theta^{-n}\omega,x_{n,k})\bigr)
\prod_{\ell=1}^{K-1}
\left(
\frac{D^\ell\Psi(j,\theta^{-n}\omega,x_{n,k})}
     {\bigl(D\Psi(j,\theta^{-n}\omega,x_{n,k})\bigr)^\ell}
\right)^{s_\ell},
\]
where
\[
\sum_{\ell=1}^{K-1}\ell s_\ell=K,
\qquad
t=\sum_{\ell=1}^{K-1}s_\ell,
\qquad
2\leq t\leq K-1.
\]

{Each quotient $D^\ell\Psi/(D\Psi)^\ell$ in this product has absolute value at most $C_2$.} Moreover, the quotient rule gives, for $1\leq\ell<K$,
\[
\left|
\frac{\mathrm{d}}{\mathrm{d}x}
\left(
\frac{D^\ell\Psi(j,\theta^{-n}\omega,x_{n,k})}
     {\bigl(D\Psi(j,\theta^{-n}\omega,x_{n,k})\bigr)^\ell}
\right)
\right|
\leq
(C_1+\ell C_2^2)
\frac{D\Psi(j,\theta^{-n}\omega,x_{n,k})}
     {J_{n,\omega}(x_{n,k})}.
\]
{The factor containing $D^t\psi_{\theta^{j-n}\omega}$ and all powers of the derivative ratios
therefore has absolute value at most $U C_2^t$, and its derivative has absolute value at most}
\[
\begin{aligned}
&U C_2^t
\left(1+t\frac{C_1}{C_2}+KC_2\right)
\frac{D\Psi(j,\theta^{-n}\omega,x_{n,k})}
     {J_{n,\omega}(x_{n,k})}
\leq 2KU C_1C_2^t.
\end{aligned}
\]
Here $t+1\leq K$, so only derivatives through order $K$ occur.

Summing these bounds with the coefficients in \eqref{lmeq4-10-3} gives $U C_3$ for the low-order factors and $2KU C_1C_3$ for their derivatives.
{Combining these two bounds with
$|D_xA^j_{n,k}|\leq C_4A^j_{n,k}$, we obtain}
{\begin{equation}\label{eq:remainder-lipschitz}
\begin{aligned}
|D(G_{n,K,\omega}-g_{n,K,\omega})(x)|
&\leq
C\sum_{j=0}^{n-1}\varepsilon^{(K-1)(n-j)}
\sum_{k=1}^{{d^n}}\frac1{J_{n,\omega}(x_{n,k})}
\leq
\frac{C{C_0}\varepsilon^{K-1}}
     {1-\varepsilon^{K-1}},
\end{aligned}
\end{equation}}
where the constant is explicitly chosen as
$
C:=U\varepsilon C_3(C_4+2KC_1).
$
The last inequality follows from Lemma~\ref{lm4.1} and the geometric series. This estimate also holds for $K=2$, with $C=0$.
Hence the remainder has a uniform Lipschitz constant.

{We next estimate the principal term $g_{n,K,\omega}$.}
Integrating the logarithmic derivative bound along the short arc joining $x$ and $y$ gives
\[
\begin{aligned}
|A^j_{n,k}(x)-A^j_{n,k}(y)|
&\leq
A^j_{n,k}(y)\bigl(e^{C_4d(x,y)}-1\bigr)
\leq
A^j_{n,k}(y)C_4e^{C_4/2}d(x,y).
\end{aligned}
\]
The chain rule and uniform expansion give
\[
0<\frac{\mathrm d}{\mathrm dx}
\Psi(j,\theta^{-n}\omega,x_{n,k})
=\frac{D\Psi(j,\theta^{-n}\omega,x_{n,k})}
       {J_{n,\omega}(x_{n,k})}
\le\varepsilon^{n-j},\qquad 0\le j\le n.
\]
Integrating a real lift along the subarc from $x$ to $y$ yields
\begin{equation}\label{eq:inverse-branch-contraction}
\begin{aligned}
d\bigl(&\Psi(j,\theta^{-n}\omega,x_{n,k}),
         \Psi(j,\theta^{-n}\omega,y_{n,k})\bigr)
\le\varepsilon^{n-j}d(x,y),
\qquad 0\le j\le n.
\end{aligned}
\end{equation}
Estimate~\eqref{eq:inverse-branch-contraction} and the $\alpha$-H\"older bound for $D^K\psi_\eta$, uniform over $\eta\in{\Omega}$, yield
\[
\begin{aligned}
&\left|
D^K\psi_{\theta^{j-n}\omega}\bigl(
            \Psi(j,\theta^{-n}\omega,x_{n,k})\bigr)
-
D^K\psi_{\theta^{j-n}\omega}\bigl(
            \Psi(j,\theta^{-n}\omega,y_{n,k})\bigr)
\right|
\leq
M\varepsilon^{\alpha(n-j)}d(x,y)^\alpha.
\end{aligned}
\]
Subtracting the products in the definition of $g^j_{x_{n,k}}$,
and using $|D^K\psi_\eta|\leq U$ for every $\eta\in{\Omega}$, gives
\[
\begin{aligned}
|g^j_{x_{n,k}}-g^j_{y_{n,k}}|
&\leq
U|A^j_{n,k}(x)-A^j_{n,k}(y)|
+
A^j_{n,k}(y)
M\varepsilon^{\alpha(n-j)}d(x,y)^\alpha
\\
&\leq
A^j_{n,k}(y)
\bigl(U C_4e^{C_4/2}+M\bigr)d(x,y)^\alpha
\leq
C_5\varepsilon^{(K-1)(n-j)}
\frac{d(x,y)^\alpha}{J_{n,\omega}(y_{n,k})},
\end{aligned}
\]
where
$
C_5:=\varepsilon\bigl(M+U C_4e^{C_4/2}\bigr).
$
We used $d(x,y)\leq d(x,y)^\alpha$ and $\varepsilon^{\alpha(n-j)}\leq1$.

Summing over $j$ and $k$, we obtain
{\begin{equation}\label{eq:principal-holder}
\begin{aligned}
|g_{n,K,\omega}(x)-g_{n,K,\omega}(y)|
&\leq
C_5d(x,y)^\alpha
\left(
\sum_{j=0}^{n-1}\varepsilon^{(K-1)(n-j)}
\right)
\sum_{k=1}^{{d^n}}\frac1{J_{n,\omega}(y_{n,k})}
\\
&\leq
\frac{C_5{C_0}\varepsilon^{K-1}}
     {1-\varepsilon^{K-1}}d(x,y)^\alpha.
\end{aligned}
\end{equation}}
By Lemma~\ref{lm4.1}, we have
\[
\sum_{k=1}^{{d^n}}\frac1{J_{n,\omega}(y_{n,k})}
=
\bigl(\mathcal L_{\theta^{-n}\omega}^{n}1\bigr)(y)
\leq{C_0}.
\]

{Combining the H\"older bound~\eqref{eq:principal-holder} for $g_{n,K,\omega}$
with the Lipschitz bound~\eqref{eq:remainder-lipschitz} for $G_{n,K,\omega}-g_{n,K,\omega}$ gives}
\[
|G_{n,K,\omega}(x)-G_{n,K,\omega}(y)|
\leq
\frac{(C+C_5){C_0}\varepsilon^{K-1}}
     {1-\varepsilon^{K-1}}d(x,y)^\alpha.
\]
The constant on the right is independent of $n,\omega,x,y$ and the choice of inverse branches.

Every pair of points at distance less than $1/2$ is contained in an arc of the type used above.
On overlapping arcs, inverse-branch labels differ only by permutations, so the sums define continuous functions on ${\mathbb{S}^1}$.
The same estimate extends to antipodal pairs by continuity. Thus the functions $G_{n,K,\omega}$ have a common $\alpha$-H\"older constant in the state variable for all $n\ge1$ and $\omega\in{\Omega}$, and are therefore equicontinuous in that variable.
This completes the proof of Lemma 4.
\end{proof}

\begin{smallprint}

\begin{ack}
This work was partially supported by NSFC (12271070, 12671219) and NSF-CQ (CSTB2023NSCQ-JQX0020).
Sandro Vaienti thanks Chongqing Normal University for the kind invitation which allowed us to work on this paper.
\end{ack}

\end{smallprint}
\end{document}